\documentclass[letterpaper]{amsart}

\usepackage[margin=1.25in]{geometry}

\usepackage{amsfonts}
\usepackage{euscript}

\usepackage{latexsym,epsfig,verbatim}
\usepackage{amsmath,amsthm,amssymb}
\usepackage{colonequals}
\usepackage{units}

\usepackage{url}
\usepackage{color}
\usepackage[colorlinks,citecolor=blue,linkcolor=red!80!black]{hyperref}
\usepackage[nameinlink]{cleveref}

\usepackage{comment}
\usepackage{amscd}
\usepackage{hyperref}
\usepackage{accents}

\usepackage{tikz}
\usetikzlibrary{matrix,arrows,decorations.pathmorphing}
\theoremstyle{plain}
\newtheorem{theorem}{Theorem}[section]
\newtheorem{lemma}[theorem]{Lemma}
\newtheorem{corollary}[theorem]{Corollary}
\newtheorem{proposition}[theorem]{Proposition}

\newtheorem*{mainthrmA}{Theorem A}
\newtheorem*{mainthrmB}{Theorem B}
\newtheorem*{mainthrmC}{Theorem C}

\newtheorem*{MainResult}{\hyperlink{thrmA}{Theorem A}}
\newtheorem*{OneSidedCTMapsSurjective}{\hyperlink{thrmB}{Theorem B}}
\newtheorem*{Surjectivity}{\hyperlink{thrmC}{Theorem C}}

\crefname{claim}{Claim}{Claims}
\newtheorem*{claim*}{Claim}

\theoremstyle{definition}

\newtheorem{defn}[theorem]{Definition}
\newtheorem{remark}[theorem]{Remark}

\newtheorem{convention}[theorem]{Convention}
\crefname{convention}{Convention}{Conventions}

\newcommand{\mb}[1]{\mathbb{#1}}
\newcommand{\mc}[1]{\mathcal{#1}}

\usepackage{enumitem}

\graphicspath{{./Pictures/}}

\begin{document}

\title[Trees, dendrites, and the Cannon-Thurston map]{\textbf{\Large Trees, dendrites, and the Cannon-Thurston map}}
\author{Elizabeth Field}
\date{\today}

\begin{abstract}
When $1\to H\to G\to Q\to 1$ is a short exact sequence of three word-hyperbolic groups, Mahan Mitra (Mj) has shown that the inclusion map from $H$ to $G$ extends continuously to a map between the Gromov boundaries of $H$ and $G$. This boundary map is known as the Cannon-Thurston map. In this context, Mitra associates to every point $z$ in the Gromov boundary of $Q$ an ``ending lamination'' on $H$ which consists of pairs of distinct points in the boundary of $H$. We prove that for each such $z$, the quotient of the Gromov boundary of $H$ by the equivalence relation generated by this ending lamination is a dendrite, that is, a tree-like topological space. This result generalizes the work of Kapovich-Lustig and Dowdall-Kapovich-Taylor, who prove that in the case where $H$ is a free group and $Q$ is a convex cocompact purely atoroidal subgroup of $\mathrm{Out}(F_N)$, one can identify the resultant quotient space with a certain $\mathbb{R}$-tree in the boundary of Culler-Vogtmann's Outer space.
\end{abstract}


\maketitle

\section{Introduction}\label{Introdection}
In \cite{CanThu07}, Cannon and Thurston showed that when $M = (S\times [0,1])/ ((x,0)\sim (\phi(x),1))$ is the mapping torus of a closed hyperbolic surface $S$ by a pseudo-Anosov homeomorphism $\phi$ of $S$, the inclusion $i: \pi_1S\rightarrow \pi_1 M$ extends to a continuous, surjective, $\pi_1(S)$-equivariant map $$\mathbb{S}^1 = \partial \mb{H}^2 = \partial \pi_1 S \xrightarrow{\partial i} \partial \pi_1 M = \partial \mb{H}^3 = \mathbb{S}^2.$$ 
Although published in 2007, this work has sparked much consideration since its circulation as a preprint in 1984. In modern terminology, if $H$ and $G$ are word-hyperbolic groups with $H\leq G$ and the inclusion map $i: H\to G$ extends to a (necessarily unique and  $H$-equivariant) continuous map between the Gromov boundaries of $H$ and $G$, $\partial i: \partial H \to \partial G$, the map $\partial i$ is called the \textit{Cannon-Thurston map}. This definition naturally extends to the more general setting of hyperbolic metric spaces. Such a map automatically exists when $H$ is an undistorted (i.e., quasiconvex) subgroup of $G$, since in that case the inclusion map is a quasi-isometric embedding. The 1984 result of Cannon and Thurston gave the first non-trivial example of the existence of such a map.

Since then, Mj (formerly Mitra) has studied the existence of Cannon-Thurston maps in settings which involve distorted subgroups of hyperbolic groups \cite{Mit98, MitTrees98, MitSurvey98}. In particular, Mitra showed in \cite{Mit98} that when 
\begin{equation}\tag{*}\label{SES}
1\to H\to G\to Q\to 1
\end{equation} is a short exact sequence of infinite word-hyperbolic groups, the Cannon-Thurston map $\partial i: \partial H\to \partial G$ exists. Since an infinite normal subgroup of infinite index in a word-hyperbolic group $G$ is not quasiconvex \cite{Ghy91}, this result gives another non-trivial example of the existence of Cannon-Thurston maps. It has been shown by Kapovich and Short \cite{KapSho96} that when $H$ is an infinite normal subgroup of a hyperbolic group $G$, the limit set of $H$ in $\partial G$ is all of $\partial G$. As this limit set is precisely the image of $\partial H$ under $\partial i$, it follows that the Cannon-Thurston map is surjective in this setting.

In \cite{Mit97}, Mitra developed a theory of ``algebraic ending laminations'' for hyperbolic group extensions to describe when points in $\partial H$ are identified under the Cannon-Thurston map $\partial i$ in the setting described above. This work provides an analog of the theory of ending laminations in the context of pseudo-Anosov homeomorphisms of surfaces developed by Thurston \cite{FLP}. To each point $z\in \partial Q$, Mitra associates an ``algebraic ending lamination'' on $H$, $\Lambda_z\subseteq \partial^2 H$, where $\partial^2H = (\partial H\times \partial H) - \mathrm{diag}$. The main result of \cite{Mit98} states that two distinct points $p,q\in \partial H$ are identified under the Cannon-Thurston map if and only if there exists some $z\in \partial Q$ for which $(p,q)$ is a leaf of the ending lamination $\Lambda_z$.

If $H$ is a torsion-free, infite-index, word-hyperbolic, normal subgroup of a word-hyperbolic group $G$, it follows from combined work of Mosher \cite{Mos96}, Paulin \cite{Pau91}, Rips-Sela \cite{RipSel94}, and Bestvina-Feighn \cite{BesFei92} that $H$ must be a free product of free groups and surface groups. For a brief explanation of this, see \cite{Mit97}. Suppose $H$ is the fundamental group of a closed hyperbolic surface $S$ and $\Gamma$ is a convex cocompact subgroup of $\mathrm{Mod}(S)$ (and hence $\Gamma$ is word-hyperbolic \cite{FarMos02}). Then, $\Gamma$ naturally gives rise to a short exact sequence $1\to H\to E_\Gamma\to \Gamma\to 1$ coming from Birman's short exact sequence for $S$. Hamenst\"{a}dt has shown that in this setting, the extension group $E_\Gamma$ is hyperbolic and the orbit map of $\Gamma$ into the curve complex of $S$ is a quasi-isometric embedding \cite{Ham05}. Since the boundary of the curve complex consists of ending laminations on $S$ \cite{Kla18}, it follows that to each point $z\in \partial \Gamma$, there is an associated ending lamination $L_z$ on the surface $S$. Mj and Rafi \cite{MjRafi18} showed that the algebraic ending lamination $\Lambda_z$ is the same as the diagonal closure of the surface lamination $L_z$. To each such ending lamination $L_z$, there is an associated dual $\mb{R}$-tree $T_z$ which can be constructed by lifting $L_z$ to $\widetilde{S}$ and collapsing each leaf and complementary component to a point. For more details, see for example \cite{Bes02, CouHilLus082}.

In the free group setting, Mitra's algebraic ending laminations for hyperbolic extensions of free groups are closely related to the theory of algebraic laminations on free groups developed by Coulbois, Hilion, and Lustig in \cite{CouHilLus082}. For any subgroup $\Gamma\leq \mathrm{Out}(F_N)$, the full preimage of $\Gamma$ under the quotient map $\mathrm{Aut}(F_N) \to \mathrm{Out}(F_N)$, also denoted by $E_\Gamma$, fits into the short exact sequence $1\to F_N\to E_\Gamma\to \Gamma\to 1$. The main result of \cite{DowTay17} states that whenever $\Gamma\leq \mathrm{Out}(F_N)$ is a convex cocompact and purely atoroidal subgroup, the extension group, $E_\Gamma$, is word-hyperbolic. In \cite{DowKapTay16}, Dowdall, Kapovich, and Taylor study the fibers of the Cannon-Thurston map $\partial i: \partial F_N\to \partial E_\Gamma$ in the case where $\Gamma\leq\mathrm{Out}(F_N)$ is convex cocompact and purely atoroidal. Since $\Gamma$ is convex cocompact, the orbit map to the free factor complex, $\mathcal{F}$, is a quasi-isometric embedding \cite{HamHen18} and hence, extends to a continuous embedding $\partial \Gamma \to \partial \mathcal{F}$. By work of Bestvina-Reynolds \cite{BesRey15} and Hamenst\"adt \cite{Ham12}, $\partial \mathcal{F}$ consists of equivalence classes of arational $F_N$-trees. Therefore, there is a class of arational $F_N$-trees, $T_z$, associated to each point $z\in \partial \Gamma$. Moreover, each such tree $T_z$ comes equipped with the ``dual lamination'' $L(T_z)$, defined by Coulbois, Hilion, and Lustig in \cite{CouHilLus082}. A key result of \cite{DowKapTay16} states that for each $z\in \partial \Gamma$, $\Lambda_z = L(T_z)$. This theorem extends the result of Kapovich and Lustig \cite{KapLus15} who prove this equality for the specific case where $\Gamma = \langle \varphi \rangle$ is the cyclic group generated by a fully irreducible, atoroidal automorphism of $F_N$. 

Given an $\mb{R}$-tree $T$, Coulbois, Hilion, and Lustig define a suitable topology on $\widehat{T} = \overline{T}\cup \partial T$, where $\overline{T}$ denotes the metric completion of $T$ and $\partial T$ is the Gromov boundary. This topology, known as the ``observers' topology'', is coarser than the Gromov topology and ensures that $\widehat{T}$ is compact. Recall that a \textit{dendrite} is a compact, connected, locally connected metrizable space which contains no simple closed curves. Coulbois, Hilion, and Lustig show that for any $\mb{R}$-tree $T$, $\widehat{T}$ equipped with the ``observers' topology'' is a dendrite, as well as a proper, Hausdorff metric space \cite{CouHilLus07}. Dendrites naturally arise from this compactification of simplicial trees, but in general can be much more complicated spaces such as certain Julia sets. Combining the result of \cite{DowKapTay16} with a general result from \cite{CouHilLus07} implies that for each $z\in \partial \Gamma$, for $\Gamma$ a convex-cocompact and purely atoroidal subgroup of $\mathrm{Out}(F_N)$, $\partial F_N / \Lambda_z$ equipped with the quotient topology is homeomorphic to $\widehat{T}_z$ equipped with the ``observers' topology''. In particular, $\partial F_N / \Lambda_z$ is homeomorphic to a dendrite. Here, $\partial F_N / \Lambda_z$ means the quotient space of $\partial F_N$ by the equivalence relation on $\partial F_N$ generated by $\Lambda_z\subseteq \partial F_n \times \partial F_n$. The main result of the present paper extends this result as follows.

\begin{MainResult}
Let $1\to H\to G \to Q \to 1$ be a short exact sequence of infinite, finitely generated, word-hyperbolic groups. For each $z\in \partial Q$, let $\Lambda_z$ denote the algebraic ending lamination on $H$ associated to $z$. Then for each $z\in \partial Q$, the space $\partial H / \Lambda_z$ is homeomorphic to a dendrite.
\end{MainResult}

We now sketch the proof of Theorem A. Let $P: \Gamma_G\to \Gamma_Q$ denote the map which is induced by the quotient map $P: G\to Q$. Let $z\in \partial Q$ be arbitrary and take any $z'\in \partial Q$ with $z'\neq z$. Consider a bi-infinite geodesic $\gamma = (z',z)\subseteq \Gamma_Q$ and define the space $X(\gamma)$ to be the subgraph of $\Gamma_G$ given by $X(\gamma) = P^{-1}(\gamma)$. We show that $X(\gamma)$ satisfies the properties of being a metric graph bundle, as defined by Mj-Sardar \cite{MjSardar12}, and that $X(\gamma)$ is hyperbolic (\Cref{X(gamma) hyperbolic}). We go on to show that $X(\gamma)$ also satisfies the properties of being a bi-infinite hyperbolic stack, as defined by Bowditch \cite{Bow13}, with fibers being copies of the Cayley graph of $H$ (\Cref{X(gamma) stack}). We then look at the semi-infinite stack $X(\gamma)^+$ which lies over the geodesic ray $\gamma^+ = [z_0,z)$, where $z_0\in (z',z)$. We denote the natural ``0-th slice'' map from $\Gamma_H\to X(\gamma)^+$ by $i_\gamma^+$, and also refer to the continuous extension of this map to $\partial i_\gamma^+: \partial H \to \partial X(\gamma)^+$ as the Cannon-Thurston map. We then show the following.

\begin{OneSidedCTMapsSurjective}
Let $1\to H\to G\to Q\to 1$ be a short exact sequence of infinite, finitely generated, word-hyperbolic groups. Let $z,z'\in \partial Q$ be distinct and let $\gamma \subseteq \Gamma_Q$ be a bi-infinite geodesic in $\Gamma_Q$ between $z$ and $z'$. Let $i_\gamma^+: \Gamma_H\to X(\gamma)^+$ be the inclusion of $\Gamma_H$ into the semi-infinite stack $X(\gamma)^+$ over $\gamma^+ = [z_0,z)$ for some $z_0\in \gamma$, and let $i_\gamma: \Gamma_H\to X(\gamma)$ be the inclusion of $\Gamma_H$ into the bi-infinite stack $X(\gamma)$ over $\gamma$. Then, 
\begin{enumerate}
    \item the Cannon-Thurston map $\partial i_\gamma^+: \partial H\to \partial X(\gamma)^+$ is surjective; and
    \item the Cannon-Thurston map $\partial i_\gamma: \partial H\to \partial X(\gamma)$ is surjective.
\end{enumerate}
\end{OneSidedCTMapsSurjective}

Using the work of Mitra from \cite{Mit97}, we show that the following holds.

\begin{Surjectivity}
Let $1\to H \to G\to Q \to 1$ be a short exact sequence of infinite, finitely generated, word-hyperbolic groups. Let $z,z'\in \partial Q$ be distinct and let $\gamma\subseteq \Gamma_Q$ be a bi-infinite geodesic between $z$ and $z'$. Let $i_\gamma^+: \Gamma_H\to X(\gamma)^+$ be the inclusion of $\Gamma_H$ into the semi-infinite stack $X(\gamma)^+$ over $\gamma^+ = [z_0,z)$ for some $z_0\in \gamma$, and let $\partial i_\gamma^+: \partial H\to \partial X(\gamma)^+$ be the Cannon-Thurston map. 

Then for any distinct $u,v\in \partial H$, we have $\partial i_\gamma^+(u) = \partial i_\gamma^+(v)$ if and only if $(u,v)$ is a leaf of the ending lamination $\Lambda_z$.
\end{Surjectivity}

To finish the proof of \hyperlink{thrmA}{Theorem A}, note that by a general result of Bowditch \cite{Bow13}, $\partial X(\gamma)^+$ is a dendrite (\Cref{BowditchDendrite}). \hyperlink{thrmC}{Theorem C} implies that the Cannon-Thurston map $\partial i_\gamma^+: \partial H\to \partial X(\gamma)^+$ quotients through to an injective map $\tau_z: \partial H / \Lambda_z \to \partial X(\gamma)^+$. Since $\partial i_\gamma^+$ is continuous, the map $\tau_z$ is also continuous. By \hyperlink{thrmB}{Theorem B}, $\tau_z$ is also surjective. Thus, $\tau_z: \partial H / \Lambda_z \to \partial X(\gamma)^+$ is a continuous bijection between two compact topological spaces, where $\partial X(\gamma)^+$ is Hausdorff. Therefore, $\tau_z$ is a homeomorphism.

In \Cref{Background}, we provide background on hyperbolic metric spaces and hyperbolic groups. The space $X(\gamma)$ is introduced in \Cref{MetricGraphBundles} and is shown to be hyperbolic. In \Cref{StacksOfSpaces}, we show that $X(\gamma)$ is a bi-infinite, hyperbolic stack and use this to prove \hyperlink{thrmB}{Theorem B}. The ending lamination $\Lambda_z$ is defined in \Cref{EndingLaminations} and several technical results are given which lead to the proof of \hyperlink{thrmC}{Theorem C}. Finally, \hyperlink{thrmA}{Theorem A} is proved in \Cref{MainResult}.

\subsection*{Acknowledgements.}
The author is very grateful to her PhD advisor Ilya Kapovich for his guidance, encouragement, feedback, and constant support. The author would also like to thank Spencer Dowdall and Sam Taylor for enlightening conversation and suggesting a generalization of an earlier version of this work, as well as Chris Leininger for his support and helpful discussions. The author would also like to thank the referee for helpful suggestions. The author gratefully acknowledges support from the NSF grant DMS-1905641 and would like to thank Hunter College for their warm hospitality during the semester this paper was written.

\section{Background}\label{Background}
In this section, we will discuss some basic definitions and facts about hyperbolic metric spaces and hyperbolic groups. For general references on hyperbolic spaces, groups, and their boundaries, see \cite{AloBra91, BriHae10, CooDelPap90, GhyHar90, Ghy91, Gro87, KapBen02}.

\subsection{Hyperbolic metric spaces.}

Let $(X,d)$ be a geodesic metric space. For any $x,y\in X$, we will denote a geodesic between $x$ and $y$ by $[x,y]_X$, or by $[x,y]$ if the space is clear. Given any three points $x,y,z\in X$, the \textit{Gromov product} of $x$ and $y$ relative to $z$ is defined to be 
\[ (x,y;X)_z := \frac{1}{2}\big( d(x,z) + d(y,z) - d(x,y)  \big). \]
If the space $X$ is clear, we will simply write $(x,y)_z$ for $(x,y;X)_z$. 

Let $\delta\geq 0$. A geodesic metric space $(X,d)$ is called \textit{$\delta$-hyperbolic} if for any $x,y,z\in X$ and any geodesics $[z,x]$ and $[z,y]$ in $X$ the following holds. Let $x'\in [z,x]$ and $y'\in [z,y]$ be any points such that $d(z,x') = d(z,y')\leq (x,y)_z$. Then, $d(x',y')\leq \delta$. Note that this property implies that for any geodesic triangle $\Delta = [x,y]\cup [y,z]\cup [z,x]$ in $X$, each side of $\Delta$ is contained in the $\delta$-neighborhood of the union of the other two sides. See \cite{AloBra91} and \cite{BriHae10} for more details and other equivalent definitions of hyperbolicity. The metric space $(X,d)$ is said to be \textit{hyperbolic} if it is $\delta$-hyperbolic for some $\delta\geq 0$. Note that in a hyperbolic metric space, the Gromov product $(x,y)_z$ measures how closely the geodesics $[z,x]$ and $[z,y]$ travel. 

A sequence of points $(x_n)_{n\in \mb{N}}\in X$ is said to \textit{converge to infinity} if for some basepoint $x\in X$, $$\liminf_{i,j\to \infty} (x_i,x_j)_x = \infty.$$ It is known that this definition is independent of basepoint. Two sequences $(x_n)$ and $(y_n)$ in $X$ which converge to infinity are said to be \textit{equivalent} if $$\liminf_{i,j\to\infty} (x_i,y_j)_x = \infty.$$ We denote the equivalence class of a sequence $(x_n)$ converging to infinity by $[(x_n)]$ and again note that this equivalence is independent of chosen basepoint. The \textit{Gromov boundary} of $X$ is defined to be 
$$\partial X := \{ [(x_n)] \mid (x_n) \text{ is a sequence converging to infinity in $X$}\}.$$ We can also represent $\partial X$ by equivalence classes of geodesic rays, where two rays represent the same point at infinity if they have bounded Hausdorff distance.

If $X$ is a proper hyperbolic metric space, then $\partial X$ is known to be compact, and so the space $\widehat{X} = X \cup \partial X$ can be considered a compactification of $X$. There is a natural topology that is carried by $\partial X$ which can be extended to a topology on $\widehat{X}$. Fix a basepoint $x\in X$ and for any $p\in \partial X$ and $r\geq 0$, define the set \begin{align*} U(p,r):= &\{q\in \partial X \mid \text{ there exist sequences } (x_n) \text{ and } (y_n) \text{ with } [(x_n)] = p \\ &\text{ and } [(y_n)] = q \text{ such that }\liminf_{i,j\to \infty} (x_i, y_j)_x\geq r\}.
\end{align*} The topology on $\partial X$ is then generated by $\{U(p,r) \mid r\geq 0\}$. To get a topology on $\widehat{X}$, we define for each $p\in \partial X$ and $r\geq 0$ the additional sets \begin{align*}
    U'(p,r):= &\{y\in X \mid \text{ for some sequence } (x_n) \text{ with } \\
    &[(x_n)] = p, \text{ we have } \liminf_{i\to\infty} (x_i, y)_x\geq r\}.
\end{align*} For each $p\in \partial X$ we put the basis of neighborhoods for $p \in \widehat{X}$ to be $\{U(p,r) \cup U'(p,r) \mid r\geq 0\}$. For each $y\in X$, we use the same neighborhood basis as in $X$. For a proper hyperbolic space, these topologies can be equivalently defined in terms of geodesic rays. Informally, two points $a,b\in\widehat{X}$ are close if geodesic rays which begin at some basepoint $x$ and end at $a$ and $b$ stay uniformly Hausdorff close for a long time. Both formulations of $\partial X$ are known to be independent of basepoint. For more details, see \cite{KapBen02}.

Let $(X, d_X)$ and $(Y, d_Y)$ be metric spaces, and let $\kappa\geq 1$ and $\epsilon\geq 0$. A map $f: X\to Y$ is said to be a \textit{$(\kappa, \epsilon)$-quasi-isometric embedding} if for all $x_1, x_2\in X$, 
$$\frac{1}{\kappa}d_X(x_1, x_2) - \epsilon \leq d_Y(f(x_1), f(x_2)) \leq \kappa d_X(x_1, x_2) + \epsilon.$$
A \textit{$(\kappa, \epsilon)$-quasigeodesic} in a metric space $(X,d)$ is the image of a $(\kappa, \epsilon)$-quasi-isometric embedding $f: I\to X$, where $I\subseteq \mathbb{R}$ is a sub-interval. The map $f$ itself is also referred to as a $(\kappa,\epsilon)$-quasigeodesic. It is known that quasigeodesics ``diverge exponentially'' in a hyperbolic metric space:

\begin{proposition}[Mitra \cite{Mit97} Proposition 2.4]\label{Mitra 2.4}
Let $(X,d)$ be a $\delta$-hyperbolic, geodesic metric space. Given $K\geq 1$, $\epsilon\geq 0$, and $\alpha\geq 0$, there exist $b > 1$, $A > 0$, and $C > 0$ such that the following holds:

If $r_1$ and $r_2$ are two $(K, \epsilon)$-quasigeodesics in $X$ with $d(r_1(0), r_2(0))\leq \alpha$ and there exists $T\geq 0$ with $d(r_1(T), r_2(T))\geq C$, then any path joining $r_1(T + t)$ to $r_2(T + t)$ and lying outside the union of the $\frac{T + t - 1}{K + \epsilon}$-balls around $r_1(0)$ and $r_2(0)$ has length greater than $Ab^t$ for all $t \geq 0$.
\end{proposition}

The following are basic facts that we will need later about hyperbolic metric spaces.

\begin{proposition}\label{OneSmallGromovProductQuasigeodesic}
Let $(X,d)$ be a $\delta$-hyperbolic metric space and let $A\geq 0$. If $x,y,z\in X$ are such that $(x,z)_y\leq A$, then $[x,y]\cup [y,z]$ is a $(1,2A)$-quasigeodesic.
\end{proposition}

\begin{proof}
Suppose that $x$, $y$, and $z$ are such that $(x,z)_y\leq A$. We need to show that for all $p\in [x,y]$ and $q\in [y,z]$, $d(p,y) + d(y,q)\leq d(p,q) + 2A$. By the triangle inequality, \begin{align*}
    (p,q)_y &= \frac{1}{2}(d(p,y) + d(q,y) - d(p,q)) \\
    &\leq \frac{1}{2}(d(p,y) + d(q,y) - (d(p,z) - d(q,z))) \\
    &= \frac{1}{2}(d(p,y) + d(y,z) - d(p,z)) = (p,z)_y.
\end{align*} Similarly, $(p,z)_y\leq (x,z)_y$. Therefore, $(p,q)_y\leq A$ by hypothesis, and so $d(p,q) + d(y,q) = d(p,q) + 2(p,q)_y \leq d(p,q) + 2A$. Hence, $[x,y]\cup [y,z]$ is a $(1,2A)$-quasigeodesic.
\end{proof}

The next proposition says that geodesic quadrilaterals in hyperbolic metric spaces must either be ``tall and thin'' or ``short and long''.

\begin{proposition}\label{ThinQuadrilateral}
Let $(X,d)$ be a $\delta$-hyperbolic metric space and let $x,y,z,w\in X$. Then, either there are points $a\in [x,y]$ and $a'\in [z,w]$ with $d(a,a')\leq 2\delta$, or there are points $b\in [x,w]$ and $b'\in [y,z]$ with $d(b,b')\leq 2\delta$.
\end{proposition}

\begin{proof}
Consider the geodesic quadrilateral with sides $[x,y]$, $[y,z]$, $[z,w]$, and $[x,w]$. Draw in the diagonal $[y,w]$ and consider the two triangles $xyw = [x,y]\cup [y,w]\cup [w,x]$ and $ywz = [y,w]\cup [w,z]\cup [z,y]$. Mark internal points $p\in [x,y]$, $q\in [x,w]$, and $r\in [y,w]$ such that $d(x,p) = d(x,q)$, $d(w,q) = d(w,r)$, and $d(y,p) = d(y,r)$. Similarly, mark internal points $q'\in [y,z]$, $p'\in [z,w]$, and $r'\in [y,w]$ such that $d(z,q') = d(z,p')$, $d(w,p') = d(w,r')$, and $d(y,q') = d(y,r')$. Note that since $X$ is $\delta$-hyperbolic, we have that $\max\{d(p,q), d(q,r), d(p,r)\}\leq \delta$ and $\max\{d(p',q'), d(q',r'), d(p',r')\}\leq \delta$. There are two cases to consider.

First, suppose that $d(y,r) \leq d(y,r')$. In this case, there exists some point $s\in [y,w]$ between $r$ and $r'$ such that $d(s,[x,w])\leq \delta$ and $d(s,[y,z])\leq \delta$. Hence, there exist $b\in [x,w]$ and $b'\in [y,z]$ such that $d(b,b')\leq d(b,s) + d(s,b') \leq 2\delta$. 

Now, suppose that $d(y,r) > d(y,r')$. In this case, there is some point $s'\in [y,w]$ between $r'$ and $r$ such that $d(s',[x,y])\leq \delta$ and $d(s',[z,w])\leq \delta$. So, there is some $a\in [x,y]$ and $a'\in [z,w]$ with $d(a,a')\leq 2\delta$. 

\end{proof}

\begin{proposition}\label{SmallGromovProductsQuasigeodesic}
Let $(X,d)$ be a $\delta$-hyperbolic metric space and let $A\geq 0$. If $x,y,z,w\in X$ are such that $(x,z)_y\leq A$, $(y,w)_z\leq A$, and $d(y,z)> 10\delta + 2A$, then $[x,y]\cup [y,z]\cup [z,w]$ is a $(1,4\delta + 4A)$-quasigeodesic.
\end{proposition}

\begin{proof}
Fix $x,y,z,w\in X$ such that $(x,z)_y\leq A$, $(y,w)_z\leq A$, and $d(y,z) > 10\delta + 2A$. We need to show that for all $p,q\in [x,y]\cup [y,z]\cup [z,w]$, the distance between $p$ and $q$ along $[x,y]\cup [y,z]\cup [z,w]$ is at most $d(p,q) + 4\delta + 4A$. This statement is certainly true if $p$ and $q$ are on the same geodesic segment, and the proof of \Cref{OneSmallGromovProductQuasigeodesic} shows that it also holds if $p$ and $q$ are on adjacent segments. So, it remains to show that if $p\in [x,y]$ and $q\in [z,w]$, then $d(p,y) + d(y,z) + d(z,q) \leq d(p,q) + 4\delta + 4A$. 

So, fix $p\in [x,y]$ and $q\in [z,w]$ and let $[p,q]$ denote the geodesic segment between $p$ and $q$. Since $d(y,z) > 10\delta + 2A$, there exists a point $r\in [y,z]$ such that $d(r,y) > 5\delta + A$ and $d(r,z) > 5\delta + A$. As geodesic quadrilaterals are $2\delta$-thin, there exists some $r'\in [y,p]\cup [p,q]\cup [q,z]$ at distance at most $2\delta$ from $r$. We claim that this point $r'\in [p,q]$. Suppose instead that $r'\in [p,y]$. Then since $(p,r)_y \leq (x,z)_y \leq A$, we have that $d(y, [p,r]) \leq A + \delta$. So, $d(z,y) \leq d(p,x) + A + \delta$. But then, \begin{align*}
    d(p,x) + A + \delta - d(x,y) &\leq d(p,r') + d(r',x) + A + \delta - d(x,r') - d(r',y) \\
    &= d(p,r') + A + \delta - d(r',y) \\
    &\leq d(p,r') + A + \delta - [d(y,p) - d(p,r')] \\
    &= 2d(p,r') + A + \delta - d(y,p) \\
    &< 4\delta + A + \delta - (5\delta + A) = 0,
\end{align*} which is a contradiction. Similarly, we cannot have that $r'\in [z,q]$ and hence our claim that $r'\in [p,q]$ must be true.

As $(p,r)_y \leq A$ and $(q,r)_z\leq A$, we have that $d(p,y) + d(y,r)\leq d(p,r) + 2A$ and $d(r,z) + d(z,q)\leq d(r,q) + 2A$. By the triangle inequality,  $d(p,r) \leq d(p,r') + 2\delta$ and $d(q,r) \leq d(q,r') + 2\delta$. Therefore, $d(p,y) + d(y,z) + d(z,q) \leq d(p,q) + 4\delta + 4A$.
\end{proof}

\begin{lemma}\label{GeodesicQuadrilateral}
Let $(X, d)$ be a $\delta$-hyperbolic metric space and let $x,y,z,w\in X$. If there exist points $a\in [x,w]$ and $b\in [y,z]$ such that $d(a,b)\leq 2\delta$, then $[x,y]\cup [y,z]\cup [z,w]$ is a $(1, 4\delta + 4d(y,z))$-quasigeodesic.
\end{lemma}

\begin{proof}
Let $x,y,z,w\in X$ be as above and consider the geodesic quadrilateral with edges $[x,y]$, $[y,z]$, $[z,w]$, and $[x,w]$. Note that both $(x,z)_y$ and $(y,w)_z$ are bounded by $d(y,z)$. So, the proof of \Cref{SmallGromovProductsQuasigeodesic} shows that if $p\in [x,y]$ and $q\in [y,z]$, then $d(p,y) + d(y,q) \leq d(p,q) + 2d(y,z)$. If $p\in [x,y]$ and $q\in [z,w]$, then there exist points $u\in [p,q]$ and $v\in [y,z]$ with $d(u,v)\leq 2\delta$. Thus, $d(p,v)\leq d(p,u) + 2\delta$ and $d(q,v) \leq d(q,u) + 2\delta$. Additionally, $d(p,y) + d(y,v) \leq d(p,v) + 2d(y,z)$ and $d(q,z) + d(z,v)\leq d(q,v) + 2d(y,z)$. Therefore, we have that \begin{align*}
    d(p,y) + d(y,z) + d(z,q) &\leq d(p,v) + d(q,v) + 4d(y,z) \\ &\leq d(p,q) + 4\delta + 4d(y,z).
\end{align*}
\end{proof}

Let $\mathcal{N}_r(U)$ denote the $r$-neighborhood around a subset $U$ of $X$. It is known that in a hyperbolic metric space, any quasigeodesic stays near the geodesic between its endpoints:
\begin{proposition}[\cite{Gro87} 7.2 A; \cite{CooDelPap90} 3.1.3; \cite{GhyHar90} 5.6, 5.11]\label{QuasigeodesicStability}
Let $(X,d)$ be a $\delta$-hyperbolic metric space and let $x,y\in X$. For any $\kappa\geq 1$ and $\epsilon\geq 0$, there exists $L = L(\delta, \kappa, \epsilon)\geq 0$ such that if $\alpha$ is a $(\kappa, \epsilon)$-quasigeodesic between $x$ and $y$, then for any geodesic $\beta = [x,y]$, we have that $\alpha \subset \mathcal{N}_L(\beta)$ and $\beta\subset \mathcal{N}_L(\alpha)$. 
\end{proposition}

\subsection{Hyperbolic groups.}

A finitely generated group $H$ is said to be \textit{word-hyperbolic} if for some, equivalently any, finite generating set of $H$, there exists $\delta\geq 0$ such that the Cayley graph of $H$ with respect to the word metric is $\delta$-hyperbolic. Let $H$ be a word-hyperbolic group and fix a finite generating set $S_H$ for $H$. We will denote the Cayley graph of $H$ with respect to $S_H$ by $\Gamma_H$ and let $d_H$, or simply $d$, denote the word-metric. Let $\partial H$ denote the Gromov boundary of $\Gamma_H$, and let $\widehat{\Gamma_H} = \partial H \cup \Gamma_H$ be the Gromov compactification of $\Gamma_H$. Then, $\widehat{\Gamma_H}$ is a compact, Hausdorff topological space. It is known that for a word-hyperbolic group, $\partial H$ is independent of choice of finite generating set.

We will now introduce some terminology and facts that we will use throughout this paper. Given a group $H$ with finite generating set $S_H$, let $\Sigma_H := S_H\cup S_H^{-1}$ denote the alphabet of $H$. A \textit{word} $w$ over the alphabet $\Sigma_H$ is an expression $s_1\cdots s_n$, where $s_i\in \Sigma_H$ and $n\geq 0$ (the case $n=0$ represents the empty word). We will denote the set of all finite words over $\Sigma_H$ by $\Sigma_H^*$, and will think of a word as the label of some (not necessarily geodesic) path in $\Gamma_H$. 

If $w\in \Sigma_H^*$ is the label of some path in $\Gamma_H$ from a vertex $a$ to $b$, then we will denote the group element $a^{-1}b\in H$ representing the word $w$ by $\overline{w}$. Given any element $h\in H$, we will denote the conjugacy class of $h$ in $H$ by $[h]_H$ (or simply by $[h]$ if the ambient group is clear). For a word $w\in \Sigma_H^*$, $|w|_H$ denotes the length of any path labeled by $w$ in $\Gamma_H$. The length of an element $h\in H$, also denoted by $|h|_H$, is defined to be the length of any geodesic from the identity $1_H$ to the vertex $h$ in $\Gamma_H$. We will drop the subscript if the group we are working in is clear.

For the remainder of this section, we assume that $H$ is a word-hyperbolic group with a fixed finite generating set $S_H$. We will also usually abbreviate $|h|_H$ by $|h|$ for $h\in H$. The following definitions generalize the notion of words in a free group being cyclically and almost cyclically reduced to the context of a general word-hyperbolic group.

\begin{defn}
Let $\kappa \geq 0$. An element $h\in H$ is said to be \textit{$\kappa$-almost conjugacy minimal} in $H$ if $|h|_H \leq |h'|_H + \kappa$ for all $h'\in [h]_H$. If $\kappa = 0$, then $h$ is said to be \textit{conjugacy minimal}. A geodesic $[a,aw]\subseteq \Gamma_H$ is said to be a \textit{$\kappa$-almost conjugacy minimal representative} if $\overline{w}\in H$ is $\kappa$-almost conjugacy minimal. If $[a,aw]\subseteq \Gamma_H$ is a $\kappa$-almost conjugacy minimal representative, then we will also refer to the word $w$ labeling this geodesic as a \textit{$\kappa$-almost conjugacy minimal representative}. If $\kappa = 0$, then $[a,aw]$ and $w$ are said to be \textit{conjugacy minimal representatives}.
\end{defn}

\begin{lemma}\label{ConjMinBoundedGromov}
Fix an element $h\in H$ and any constant $\kappa \geq 0$. If $h$ is $\kappa$-almost conjugacy minimal, then $(1,hh)_h \leq \frac{\kappa + \delta}{2}$.
\end{lemma}

\begin{proof}
Let $h\in H$ be $\kappa$-almost conjugacy minimal and suppose the geodesic $[1,h]\subseteq \Gamma_H$ is labeled by $\alpha h' \beta$, where $|\alpha| = |\beta| = (1,hh)_h$ and $\beta\alpha = s$, with $|s| \leq \delta$. Then $h =_H \alpha h' s \alpha^{-1}$, and so $h's\in [h]_H$. As $h$ is $\kappa$-almost conjugacy minimal, we have that $|h| \leq |h's| + \kappa \leq |h'| + \delta + \kappa$. Finally, $|h| = |\alpha| + |h'| + |\beta| = 2(1,hh)_h + |h'|$, and so we have that $(1,hh)_h\leq \frac{\kappa + \delta}{2}$.
\end{proof}

\begin{lemma}\label{GeneralWordAsFHR}
There exists a constant $C\geq 0$ such that for any element $h\in H$, the following holds. Suppose that $u,c\in H$ are such that $h = c^{-1}uc$, where $u\in [h]$ is conjugacy minimal and $|c|$ is the smallest element conjugating $h$ to any conjugacy minimal representative. Then, the path $[c,1]\cup [1,u] \cup [u,uc]\subseteq \Gamma_H$ is a $(1,C)$-quasigeodesic.
\end{lemma}

\begin{proof}
Let $h, u,c\in H$ be as in the hypothesis above and consider the quadrilateral in $\Gamma_H$ with vertices $1$, $c$, $u$, $uc$, and edges $[1,u]$ labeled by $u$, $[c,uc]$ labeled by $h$, and $[1,c]$ and $[u,uc]$ both labeled by $c$. We want to show the path $\gamma = [c,1]\cup [1,u]\cup [u,uc]$ is a $(1,C)$-quasigeodesic, for some constant $C\geq 0$. 

Let $p\in [1,c]$ and $q\in [1,u]$ be such that $d(1,p) = d(1,q) = (c,u)_1$. As $d(p,q)\leq \delta$, we must have that $d(1,p)$ is also at most $\delta$. Otherwise, $q^{-1}c$ would be a shorter word conjugating $h$ to a cyclic conjugate of $u$, contradicting the minimality of $|c|$. Similarly, we have that $(1,uc)_u\leq \delta$. If $|u| > 12\delta$, then by \Cref{SmallGromovProductsQuasigeodesic}, $\gamma$ is a $(1,8\delta)$-quasigeodesic. 

If $|u| \leq 12 \delta$, then since $H$ is finitely generated, there are only finitely many possibilities for such $u$. Hence, there are only finitely many cases to consider and the result holds by taking $C$ to be, for instance, the length of the longest path $\gamma$ that we get in this setting.
\end{proof}

\begin{corollary}\label{ShortElementToConjMin}
For any $\kappa \geq 0$, there exists a constant $M>0$ such that if $h\in H$ is $\kappa$-almost conjugacy minimal, then there is an element $c\in H$ with $|c|\leq M$ and a conjugacy minimal element $u\in [h]$ such that $h = c^{-1}uc$.
\end{corollary}

\begin{proof}
Let $c\in H$ be a shortest length element conjugating $h$ to any conjugacy minimal element in $[h]$. By \Cref{GeneralWordAsFHR}, there exists some constant $C> 0$ such that $[c,1]\cup [1,u]\cup [u,uc]$ is a $(1,C)$-quasigeodesic. So, $2|c| + |u| \leq |h| + C$. Since $h$ is $\kappa$-almost conjugacy minimal, we have that $|h| \leq |u| + \kappa$. Thus, $|c| \leq \frac{C + \kappa}{2}$.
\end{proof}

\begin{lemma}
For any $\kappa\geq 0$ there exists a constant $A\geq 0$ such that if $h\in H$ satisfies $(1,hh)_h\leq A$, then $h$ is $\kappa$-almost conjugacy minimal.
\end{lemma}

\begin{proof}
Let $h\in H$ be such that $(1,hh)_h \leq A$. Then by \Cref{OneSmallGromovProductQuasigeodesic}, the path $[1,h]\cup [h,hh]$ is a $(1,2A)$-quasigeodesic. Additionally, all subpaths of $[1,h] \cup [h,hh]$ are $(1,2A)$-quasigeodesics. In particular, any (non-reduced) edge-path representing a cyclic conjugate of $h$ is a $(1,2A)$-quasigeodesic. Choose a cyclic conjugate, $h'$ of $h$ such that $ch'c^{-1} = u$, where $u\in [h]_H$ is conjugacy minimal and $|c|$ is smallest.

Consider the points $1$, $c$, $u$, and $uc$; geodesics $[1,c]$, $[1,u]$, and $[u,uc]$; and the $(1, 2A)$-quasigeodesic path between $c$ and $uc$, call it $\gamma'$, labeled by the (non-reduced) word $h'$. 
By \Cref{GeneralWordAsFHR}, there is some constant $C$ for which $\gamma = [c,1]\cup [1,u]\cup [u,uc]$ is a $(1,C)$-quasigeodesic. As $\gamma'$ and $\gamma$ are quasigeodesics sharing the same endpoints, \Cref{QuasigeodesicStability} implies that  $\gamma'$ and $\gamma$ live in a $D$-neighborhood of each other for some constant $D\geq 0$ depending only on the quasi-isometry constants and $\delta$. 

We will now show that $|c|$ is bounded. If $|u| \leq 12\delta$, then there are only finitely many cases to check and we can take maximum length we get in these cases. So, suppose that $|u| > 12\delta$. Note that the distance between any point on $[1,c]$ must be at least $|u|$ from a point on $[u,uc]$ as otherwise we would get a contradiction with $u$ being conjugacy minimal. So by \Cref{ThinQuadrilateral}, there must exist points $x\in [1,u]$ and $x'\in [c,uc]$ such that $d(x,x')\leq 2\delta$. Let $x_0$ denote the point along $\gamma'$ where the two paths labeled by $h$ meet. As the triangle with vertices $c$, $x_0$, and $uc$ is $\delta$-thin, there must exist a point $x''\in [c,x_0] \cup [x_0, uc] = \omega'$ such that $d(x', x'')\leq \delta$. Therefore $d(x,x'')\leq 3\delta$. Now, consider the word $c'$ which labels the path from $x$ to $x''$ and note that $c'$ conjugates a cyclic conjugate of $u$ to a cyclic conjugate of $h'$. Therefore, by the minimality of $c$, we must have in this case that $|c| \leq |c'|\leq 3\delta$. 

We now want to show that the distance between $x_0$ and $[1,u]$ is bounded. Consider the point $y_0\in \gamma$ which is closest to $x_0$. Without loss of generality, we may assume that either $y_0\in [1,u]$ or $y_0\in [c,1]$. If $y_0\in [1,u]$, then $d(x_0, [1,u])\leq M$. If $y_0\in [c,1]$, then $d(x_0, [1,u])\leq M + |c|$. As $|c|$ is bounded by some constant, we have that the distance between $x_0$ and $[1,u]$ is also bounded by some constant. Therefore, $h$ is $\kappa$-almost conjugacy minimal for some $\kappa\geq 0$ independent of $h$. 
\end{proof}

For the purpose of this paper, if $X$ is a graph, then we will assume that any quasi-isometry or quasi-isometric embedding takes vertices to vertices and edges to edge-paths. The following lemma follows from \Cref{QuasigeodesicStability}.

\begin{lemma}\label{FHRQGeo}
Let $K\geq 1$ and $C\geq 0$. Then, for any $\kappa\geq 0$, there exists $\kappa'\geq 0$ such that if $w\in\Sigma_H^*$ is a $\kappa$-almost conjugacy minimal representative and $\psi: \Gamma_H\to \Gamma_H$ is any $(K, C)$-quasi-isometry, then $\psi(w)$ is a $\kappa'$-almost conjugacy minimal representative. 
\end{lemma}

\section{Metric Graph Bundles}\label{MetricGraphBundles}

In \cite{BesFei92}, Bestvina and Feighn explored the question of when a space which results from the combination of Gromov-hyperbolic spaces will itself be hyperbolic. They introduced the notion of a graph of spaces and provided a ``flaring'' condition which gives a sufficient condition for the hyperbolicity of a graph of hyperbolic spaces. Mj and Sardar generalized this work in \cite{MjSardar12} where they introduced the notion of a metric graph bundle and defined the following flaring condition. 

Let $X$ and $B$ be connected graphs, each equipped with the path metric where each edge has length 1, and let $p: X\rightarrow B$ be a simplicial surjection. For the purpose of this paper, we will consider $\mb{N} = \mb{Z}_{\geq 0}$.

\begin{defn}
$X$ is said to be a \textit{metric graph bundle} over $B$ if there exists a function $f: \mb{N}\rightarrow \mb{N}$ such that:
\begin{enumerate}[label=(B\arabic*)]
\item\label{(B1)} For each vertex $b\in V(B)$, the fiber $F_b:= p^{-1}(b)$ is a connected subgraph of $X$; and for all vertices $u,v\in V(F_b)$, the induced path metric $d_b$ on $F_b$ satisfies $d_b(u,v)\leq f(d_X(u,v))$.
\item\label{(B2)} If $b_1,b_2\in V(B)$ are any two adjacent vertices and if $x_1\in V(F_{b_1})$ is any vertex, then there is some vertex $x_2\in V(F_{b_2})$ adjacent to $x_1$ in $X$.
\end{enumerate}
\end{defn}

\begin{remark}\label{MetricGraphBundle}
Note that if $p:X\to B$ is a metric graph bundle and $W\subseteq B$ is any connected subgraph, then $p: p^{-1}(W) \to W$ is again a metric graph bundle. 
\end{remark}

Given any metric graph bundle $p: X\rightarrow B$ and a connected, closed interval $I\subseteq \mb{R}$, a \textit{$(k,k)$-quasi-isometric lift} of a geodesic $\gamma: I\rightarrow B$ is any $(k,k)$-quasigeodesic $\widetilde{\gamma}: I\rightarrow X$ for which $p(\widetilde{\gamma}(n)) = \gamma(n)$ for all $n\in I\cap \mb{Z}$. 

\begin{defn}\label{def: flaring condition}
The metric graph bundle $p: X\rightarrow B$ is said to satisfy the \textit{flaring condition} if for all $k\geq 1$, there exists $\lambda_k> 1$ and $n_k, M_k\in \mb{N}$ such that the following holds: If $\gamma: [-n_k, n_k]\rightarrow B$ is any geodesic and $\widetilde{\gamma_1}$ and $\widetilde{\gamma_2}$ are any two $(k,k)$-quasi-isometric lifts of $\gamma$ in $X$ which satisfy $d_{\gamma(0)}(\widetilde{\gamma_1}(0), \widetilde{\gamma_2}(0))\geq M_k$, then we have $$\lambda_k \cdot d_{\gamma(0)}(\widetilde{\gamma_1}(0), \widetilde{\gamma_2}(0)) \leq \max\{ d_{\gamma(-n_k)}(\widetilde{\gamma_1}(-n_k), \widetilde{\gamma_2}(-n_k)), d_{\gamma(n_k)}(\widetilde{\gamma_1}(n_k), \widetilde{\gamma_2}(n_k))\}.$$
\end{defn}

The following are two theorems of Mj and Sardar which we will use later. The first is their combination theorem for metric graph bundles, which generalizes the combination theorem of Bestvina-Feighn \cite{BesFei92}. The second shows that flaring is a necessary condition for the hyperbolicity of a metric graph bundle.

\begin{theorem}[Mj-Sardar \cite{MjSardar12}]\label{theorem: Mj-Sardar Theorem 4.3}
Suppose that $p: X\rightarrow B$ is a metric graph bundle which satisfies:
\begin{enumerate}
\item $B$ is a $\delta$-hyperbolic metric space; 
\item for each $b\in V(B)$, the fiber $F_b$ is $\delta$-hyperbolic with respect to $d_b$, the path metric induced by $X$;
\item for each $b\in V(B)$, the set of barycenters of ideal triangles in $F_b$ is $D$-dense; and
\item the flaring condition is satisfied.
\end{enumerate}
Then, $X$ is a hyperbolic metric space.
\end{theorem}

\begin{theorem}[Mj-Sardar \cite{MjSardar12}]\label{proposition: Mj-Sardar Proposition 5.8}
Suppose that $p: X\rightarrow B$ is a metric graph bundle which satisfies:
\begin{enumerate}
\item $X$ is $\delta$-hyperbolic; and
\item for each $b\in V(B)$, the fiber $F_b$ is $\delta$-hyperbolic with respect to $d_b$, the path metric induced by $X$.
\end{enumerate}
Then, the metric bundle satisfies the flaring condition.
\end{theorem}

Throughout the remainder of this paper, we will use the following conventions:

\begin{convention}\label{ShortExactSequenceConvention}
For the remainder of the paper, unless otherwise specified, let $1\rightarrow H\overset{i}{\rightarrow} G \overset{P}{\rightarrow} Q \rightarrow 1$ be a short exact sequence of three infinite, word-hyperbolic groups. Fix finite, symmetric generating sets $S_H$, $S_G$, and $S_Q$ for $H$, $G$, and $Q$, respectively, so that $i(S_H)\subseteq S_G$ and $S_Q:= P(S_G)$. Let $\Gamma_H$, $\Gamma_G$, and $\Gamma_Q$ denote the Cayley graphs with respect to these generating sets. Let $P: \Gamma_G\to \Gamma_Q$ also denote the map on the Cayley graphs induced by $P: G\to Q$ which is given as follows. If $v\in \Gamma_G$ is a vertex labeled by the element $g\in G$, then $v$ will get sent to the vertex in $\Gamma_Q$ labeled by the element $P(g)\in Q$. Suppose $e = [g_1,g_2]\in\Gamma_G$ is an edge between adjacent vertices $g_1, g_2\in \Gamma_G$. If $g_1$ and $g_2$ are in the same coset of $H$ in $G$, then $e$ will get collapsed to the vertex $P(g_1) = P(g_2)$ in $\Gamma_Q$. Otherwise, $e$ will get mapped to the edge between $P(g_1)$ and $P(g_2)$ in $\Gamma_Q$ labeled by $P(g_1^{-1}g_2)\in S_Q$.
\end{convention}

\begin{convention}\label{StackConvention}
Suppose $\gamma = (z',z)$ is a bi-infinite geodesic in $\Gamma_Q$ between $z',z\in \partial Q$ with $z'\neq z$; and let $z_0\in V(\gamma)$ be a vertex of $\gamma$ which minimizes $d_Q(1, \gamma)$. Label the sequence of vertices in order along the portion of $\gamma$ from $z_0$ to $z$ by $z_0, z_1, z_2, \ldots$; and similarly, label the sequence of vertices in order along the portion of $\gamma$ from $z_0$ to $z'$ by $z_0, z_{-1}, z_{-2}, \ldots$. Let $\gamma^+ = [z_0,z)$ and $\gamma^- = (z',z_0]$.
\end{convention}

\begin{defn}
The \textit{subgraph of $\Gamma_G$ corresponding to $\gamma$} is \[X(\gamma):= P^{-1}(\gamma). \]
\end{defn}

Note that we can think of $X(\gamma)$ as the subgraph of $\Gamma_G$ with vertical fibers that are copies of $\Gamma_H$ corresponding to the cosets $g_iH$, where $g_i\in P^{-1}(z_i)$ for each $z_i\in V(\gamma)$. Since $S_Q = P(S_G)$, there are edges between adjacent cosets $g_iH$ and $g_{i+1}H$ between any vertex $g_ih$ and the vertex $g_ihP^{-1}([z_i,z_{i+1}])$, where $[z_i,z_{i+1}]$ is the edge in $\gamma$ between $z_i$ and $z_{i+1}$. Let $P_\gamma: X(\gamma)\rightarrow \gamma$ denote the restriction of $P$ to $X(\gamma)$.

Mj and Sardar showed in \cite{MjSardar12} that $P: \Gamma_G\to \Gamma_Q$ is a metric graph bundle. The same reasoning shows that the restricted map $P_\gamma$ is a metric graph bundle as well. We include the argument below for completeness.

\begin{proposition}
Given $P: \Gamma_G\to \Gamma_Q$ as in \Cref{ShortExactSequenceConvention}, the map $P: \Gamma_G\to \Gamma_Q$ and the restricted map $P_\gamma: X(\gamma)\to \gamma$ are metric graph bundles.
\end{proposition}

\begin{proof}
For each vertex $q\in V(\Gamma_Q)$, $P^{-1}(q) = F_q$ is a copy of $\Gamma_H$, and so the induced path metric $d_q$ is equal to $d_H$ for all $q$. Hence, condition (B1) is satisfied by the function $f(n):= \max\{d_H(1,g) \mid d_G(1,g)\leq n \}$. Now, suppose $q_1,q_2\in \Gamma_Q$ are adjacent vertices where $P(g_1H) = q_1$ and $P(g_2H) = q_2$. Since $P$ maps edges between distinct cosets of $\Gamma_H$ in $\Gamma_G$ isometrically onto edges in $\Gamma_Q$, there exist some $h_1, h_2\in H$ such that $g_1h_1$ and $g_2h_2$ are adjacent in $\Gamma_G$. Therefore $s = (g_1h_1)^{-1}g_2h_2\in S_G$. Hence, for all $x_1 = g_1h\in V(F_{q_1})$, $x_1$ is adjacent to $x_1 s = g_1hs = (g_1hh_1^{-1}g_1^{-1})g_2h_2$ in $\Gamma_G$. This element is contained in the coset $g_2H = F_{q_2}$ since $H$ is normal in $G$, and so condition \ref{(B2)} is satisfied. By \Cref{MetricGraphBundle}, $P_\gamma: X(\gamma)\to \gamma$ is also a metric graph bundle.
\end{proof}

Condition \ref{(B2)} says that if we choose any lift $g_0$ of $z_0$, there exists $g_1\in P^{-1}(z_1)$ such that $d_{X(\gamma)}(g_0, g_1) = 1$. Continuing in this fashion, we get a lift $\sigma: \gamma \rightarrow X(\gamma)$, where $\sigma(z_i) = g_i$, such that $d_{X(\gamma)}(g_i, g_{i+1}) = 1$ for all $i$. By the triangle inequality and the fact that $\gamma$ is a geodesic in $\Gamma_Q$, we have that $d_{X(\gamma)}(g_i, g_j)\leq d_Q(Pg_i, Pg_j) = d_Q(z_i, z_j)$. But, as every path in $\Gamma_G$ projects to a path in $\Gamma_Q$ of no greater length and as $X(\gamma)\subseteq \Gamma_G$, we also have that $d_Q(Pa,Pb)\leq d_G(a,b)\leq d_{X(\gamma)}(a,b)$. Hence, for all $g_i,g_j\in \sigma(\gamma)$, $d_{X(\gamma)}(g_i, g_j) = d_G(g_i, g_j) = d_Q(z_i, z_j)$. 

\begin{proposition}\label{X(gamma) hyperbolic}
The space $X(\gamma)$ is hyperbolic.
\end{proposition}

\begin{proof}
By \Cref{proposition: Mj-Sardar Proposition 5.8}, $\Gamma_G$ satisfies the flaring condition since $\Gamma_G$ and $\Gamma_H$ are both hyperbolic and for each $q\in \Gamma_Q$, $F_q:= p^{-1}(q)$ is a copy of $\Gamma_H$. Suppose that $\sigma$ is a $(K,C)$-quasi-isometric lift of $\gamma$ to $X(\gamma)$. Note that for all $a,b\in \gamma$, $d_Q(a,b) = d_Q(P\cdot \sigma(a), P\cdot \sigma(b)) \leq d_G(\sigma(a), \sigma(b))$. Also, since $X(\gamma)\subseteq \Gamma_G$, $d_G(\sigma(a), \sigma(b)) \leq d_{X(\gamma)}(\sigma(a), \sigma(b))$. So, any quasi-isometric lift of a portion of $\gamma$ to $X(\gamma)$ is also a quasi-isometric lift when considered as a path in $\Gamma_G$. Thus, we have that $X(\gamma)$ satisfies the flaring condition. Additionally, the barycenters of ideal triangles in $\Gamma_H$ are dense since the $H$-orbit of the barycenter of any ideal triangle in $\Gamma_H$ is dense in $\Gamma_H$. Therefore by \Cref{theorem: Mj-Sardar Theorem 4.3}, we have that $X(\gamma)$ is hyperbolic.
\end{proof}

\section{Stacks of Spaces}\label{StacksOfSpaces}

In \cite{Bow13}, Bowditch defines the notion of a stack of spaces. We will show that the bundle $X(\gamma)$ described above can be thought of as a hyperbolic stack of spaces. 

\begin{defn}\label{definition: straight}
Let $(X, d_X)$ and $(Y, d_Y)$ be path-metric spaces. A map $f: X\rightarrow Y$ is said to be \textit{straight} if there exist functions $F_1, F_2: [0,\infty)\rightarrow [0,\infty)$ such that for all $x,x'\in X$, $F_1(d_X(x,x'))\leq d_Y(f(x),f(x')) \leq F_2(d_X(x,x'))$, where $F_1(t)\rightarrow \infty$ as $t\rightarrow \infty$. If $X\subseteq Y$, we say that $X$ is a \textit{straight} subspace if the inclusion map $i: X\rightarrow Y$ is a straight map with respect to the induced path metric on $X$.
\end{defn}

\begin{defn}
Let $(\mc{X}, \rho)$ be a geodesic space, and let $((X_i, \rho_i))_{i\in\mb{Z}}$ be a sequence of geodesic subspaces, $X_i\subseteq \mc{X}$, called the \textit{sheets} of $\mc{X}$ with uniform quasi-isometries $f_i: X_i\rightarrow X_{i+1}$. The space $(\mc{X}, \rho)$ is said to be a \textit{bi-infinite hyperbolic stack} if it satisfies the conditions (S1)-(S6) stated below.

\begin{enumerate}[label=(S\arabic*)]
\item\label{(S1)} Each of the spaces $(X_i, \rho_i)$ are uniformly straight in $\mc{X}$, and $\rho(X_i, X_j)$ is bounded away from 0 for $i\neq j$.  
\item\label{(S2)} For all $i,j\in \mb{Z}$, $\rho(X_i, X_j)$ is bounded below by an increasing linear function of $|i - j|$.
\item\label{(S3)} For all $i\in \mb{Z}$, haus$(X_i, X_{i+1})$ is bounded above.
\item\label{(S4)} The spaces $(X_i, \rho_i)$ are uniformly hyperbolic geodesic spaces.
\item\label{(S5)} The space $(\mc{X}, \rho)$ is hyperbolic.
\item\label{(S6)} The union $\bigcup_{i\in\mb{Z}} X_i$ is quasidense in $\mc{X}$.
\end{enumerate}
\end{defn}

Given a bi-infinite stack $\mc{X}$, denote by $\mc{X}^+$ and $\mc{X}^-$ the subsets of $\mc{X}$ which consist of the sheets $(X_i)_{i\in\mb{N}}$ and $(X_i)_{i\in -\mb{N}}$, respectively. Here, $\mathbb{N} = \mathbb{Z}_{\geq 0}$ and $-\mathbb{N} = \mathbb{Z}_{\leq 0}$. We will refer to $\mc{X}^+$ and $\mc{X}^-$ as \textit{semi-infinite} stacks.

\subsection{General background on stacks} 

We first give some general background on stacks of spaces which we will later apply to the space $X(\gamma)$. Bowditch proves the following about stacks of hyperbolic spaces indexed by any subset $I\subseteq \mb{Z}$ of consecutive integers.

\begin{proposition}[Bowditch \cite{Bow13} Proposition 2.1.7]
Suppose $\mc{X}$ is a bi-infinite stack with uniformly hyperbolic sheets $(X_i)_{i\in \mb{Z}}$. If $\mc{X}$ is hyperbolic, then so is $\mc{X}(I)$, where $I\subseteq \mb{Z}$ is any set of consecutive integers. In particular, the \textit{semi-infinite} stacks $\mc{X}^+$ and $\mc{X}^-$ are hyperbolic whenever $\mc{X}$ is hyperbolic. 
\end{proposition}

Given a (bi-infinite) stack $\mc{X}$, Bowditch defines an \textit{$r$-chain}, $(x_i)_{i\in\mc{I}}$, to be a sequence of points, $x_i\in X_i$, such that $\rho(x_i, x_{i+1})\leq r$ for all $i\in \mc{I}$. A \textit{bi-infinite}, \textit{positive}, and \textit{negative} $r$-chain is defined to be an $r$-chain indexed by $\mb{Z}$, $\mb{N}$, and $-\mb{N}$, respectively. Bowditch notes that each $r$-chain interpolates a quasigeodesic in $\mc{X}$. If $\mc{X}$ is a hyperbolic stack, it comes equipped with its Gromov boundary, $\partial \mc{X}$. Thus when $\mc{X}$ is a proper, hyperbolic stack, each positive and negative chain determines a point of $\partial \mc{X}$. In this setting, there is a fixed $r_0$ depending on the hyperbolicity constant of $\mc{X}$ for which each point in $\mc{X}$ is contained in some $r_0$-chain. Bowditch defines $\partial^+ \mc{X}$ (respectively $\partial^- \mc{X}$) to be those subsets of $\partial \mc{X}$ which are determined by positive (respectively negative) $r_0$-chains. Note that the positive chains in $\mc{X}^+$ are exactly the positive chains in $\mc{X}$, and the negative chains in $\mc{X}^-$ are exactly the negative chains in $\mc{X}$. Furthermore, two chains determine the same point in $\partial\mc{X}^+$ or $\partial\mc{X}^-$ if and only if those two chains determine the same point in $\partial\mc{X}$. Hence on the level of sets, we can identify $\partial^+ \mc{X}^+$ with $\partial^+ \mc{X}$ and $\partial^-\mc{X}^-$ with $\partial^-\mc{X}$.

Each of the sheets $X_i$ are quasi-isometric to one another, and so we get a homeomorphism from $\partial X_i$ to $\partial X_j$, for all $i,j\in\mb{Z}$. We will let $\partial X_0$ denote this space which is homeomorphic to $\partial X_i$ for all $i\in \mb{Z}$. The notion of the Cannon-Thurston map, as defined earlier between the boundaries of hyperbolic groups, can be extended in the natural way to be defined between the boundaries of hyperbolic spaces. Bowditch proves the following statements about the Cannon-Thurston maps in this setting of stacks of spaces.

\begin{proposition}[Bowditch \cite{Bow13} see 2.3.2 and 2.3.3]\label{Bowditch Prop 2.3.2} \label{Bowditch Prop 2.3.3} Let $\mc{X}$ be a bi-infinite hyperbolic stack, let $\mc{X}^+$ and $\mc{X}^-$ be semi-infinite proper hyperbolic stacks, and let $\omega$, $\omega^+$, and $\omega^-$ denote the inclusions of $X_0$ into $\mc{X}$, $\mc{X}^+$, and $\mc{X}^-$, respectively. Then, 
\begin{enumerate}
\item The following continuous Cannon-Thurston maps exist: $\partial\omega: \partial X_0 \to \partial \mc{X}$, $\partial\omega^+: \partial X_0 \to \partial\mc{X}^+$, and $\partial\omega^-: \partial X_0 \to \partial \mc{X}^-$;
\item $\partial \mc{X} = \partial^+ \mc{X} \cup \partial^- \mc{X}\cup \partial\omega(\partial X_0)$; and
\item $\partial \mc{X}^+ = \partial^+ \mc{X} \cup \partial\omega^+(\partial X_0)$ and $\partial \mc{X}^- = \partial^- \mc{X} \cup \partial\omega^-(\partial X_0)$.
\end{enumerate}
\end{proposition}

Given the Cannon-Thurston maps $\partial \omega$ and $\partial\omega^\pm$, denote by $\widehat{\omega}$ and $\widehat{\omega}^\pm$ the continuous extensions of the inclusion maps. Bowditch defines the maps $\partial\tau^\pm: \partial \mc{X}^\pm \rightarrow \partial\mc{X}$ which extend to continuous maps $\widehat{\tau}^\pm: \widehat{\mc{X}^\pm}\to \widehat{\mc{X}}$ such that $\widehat{\omega} = \widehat{\tau}^\pm \circ \widehat{\omega}^\pm$. For $y\in \partial^+\mc{X}^+ = \partial^+ \mc{X}$, the map $\partial\tau^+$ is given by $\partial\tau^+(y) = y$; and for $a\in \partial X_0$, we have that $\partial\tau^+\circ \partial\omega^+(a) = \partial\omega(a)$. The map $\partial\tau^-$ is defined similarly. Bowditch proves that $\partial\tau^\pm$ are continuous maps. Using this structure, Bowditch shows the following.

\begin{lemma}[Bowditch \cite{Bow13} see 2.3.5, 2.3.6, 2.3.7, and 2.3.9] \label{Bowditch Lemma 2.3.5} \label{Bowditch Prop 2.3.6} \label{Bowditch Prop 2.3.7} \label{Bowditch Prop 2.3.9} Let $\mc{X}$ be a bi-infinite, proper, hyperbolic stack. 
\begin{enumerate}
\item Suppose $a\in \partial X_0$ and $y\in \partial^+ \mc{X}$. Then, $\partial\omega(a) = y$ if and only if there is a sequence $(\underline{x}^n)_{n\in\mb{N}}$ of positive chains, $\underline{x}^n = (x_i^n)_{i\in\mb{N}}$, each converging to $y$, and with $x_0^n$ converging to $a\in \partial X_0$.
\item Given $a\in \partial X_0$ and $y\in \partial^\pm \mc{X}$, we have $\partial\omega^\pm(a) = y$ if and only if $\partial\omega(a) = y$.
\item Suppose $a,b\in\partial X_0$ are distinct. If $\partial\omega^+(a) = \partial\omega^+(b) = y$, then $y\in \partial^+ \mc{X}$; and if $\partial\omega^-(a) = \partial\omega^-(b) = y$, then $y\in \partial^- \mc{X}$.
\item If $a,b\in \partial X_0$ and $\partial\omega(a) = \partial\omega(b)$, then either $\partial\omega^+(a) = \partial\omega^+(b)$ or $\partial\omega^-(a) = \partial\omega^-(b)$.
\end{enumerate}

\end{lemma}

\subsection{Application of stacks} We now apply this work of Bowditch to our setting of hyperbolic group extensions. Let $\gamma$ be as in \Cref{StackConvention}, and recall that $P: \Gamma_G\to \Gamma_Q$ is the projection map and $X(\gamma):= P^{-1}(\gamma)$. 

\begin{proposition}\label{X(gamma) stack}
The space $X(\gamma)$ with the induced path metric $d_{X(\gamma)}$ from $\Gamma_G$ is a hyperbolic stack.
\end{proposition}

\begin{proof}
We need to show that $X(\gamma)$ satisfies conditions \ref{(S1)}-\ref{(S6)}. For each vertex $z_i\in \gamma$, choose some $g_i\in G$ such that $P(g_i) = z_i$. For each $i\in \mb{Z}$, the sheet $X_i$ of $X(\gamma)$ is the copy of $\Gamma_H$ which corresponds to the coset $g_i\Gamma_H$ of $H$ in $G$. Since $X_i$ and $X_j$ represent different cosets of $\Gamma_H$ in $\Gamma_G$ for $i\neq j$, we have that $d_G(X_i, X_j)\leq d_{X(\gamma)}(X_i, X_j)$ is bounded away from 0 for $i\neq j$. Now, for all $i\in\mb{Z}$, let $\beta_i(n) := \max\{ d_{X_i}(a,b) \mid d_{X(\gamma)}(a,b)\leq n\}$. Then, $\beta_i^{-1}(d_{X_i}(a,b))\leq d_{X(\gamma)}(a,b) \leq d_{X_i}(a,b)$, and so condition \ref{(S1)} is satisfied.

We see that condition \ref{(S2)} is satisfied since $d_{X(\gamma)}(X_i, X_j)\geq d_Q(z_i, z_j) = |i - j|$. Similarly, we have that the Hausdorff distance between $X_i$ and $X_{i+1}$ in $X(\gamma)$ is at most 2, and so condition \ref{(S3)} is satisfied. As each $X_i$ is a copy of $\Gamma_H$ which is $\delta$-hyperbolic, we have that \ref{(S4)} holds. Additionally, $\bigcup_{i\in\mb{Z}} X_i$ is in the 1-neighborhood of $X(\gamma)$, and so \ref{(S6)} is satisfied. Finally, we have by \Cref{X(gamma) hyperbolic} that condition \ref{(S5)} is satisfied. Therefore, we have that $X(\gamma)$ is a bi-infinite hyperbolic stack.
\end{proof}

Recall, as in \Cref{StackConvention}, that $z_0$ denotes a point on $\gamma$ closest to the identity in $\Gamma_Q$, the vertices along $\gamma$ between $z_0$ and $z$ are labeled by $z_1, z_2, \ldots$, and the vertices along $\gamma$ between $z_0$ and $z'$ are labeled by $z_{-1}, z_{-2}, \ldots$. Then, for all $x_i\in X_i$, $Px_i = z_i$. Since $X(\gamma)$ satisfies property \ref{(B2)} of being a metric graph bundle, every vertex in $X(\gamma)$ is contained in some 1-chain. 
Let $y\in \partial X(\gamma)$, and let $y_n\in X(\gamma)$ be a sequence of vertices in $X(\gamma)$ which converge to $y$. As every vertex in $X(\gamma)$ is contained in some 1-chain, for each $n\in \mathbb{N}$ we can construct a 1-chain $\underline{x}^n = (x_i^n)_{i=0}^{m_n}$ in $X(\gamma)$ with terminal point $x_{m_n}^n := y_n$ as follows. Without loss of generality, assume that $y_n\in X(\gamma)^+$. Then, there exists some $m_n\in \mathbb{N}$ and some $h\in H$ such that $y_n = g_{m_n}h\in X_{m_n} = g_{m_n}\Gamma_H$, where $g_{i} = \sigma(z_{i})$. Set $x_{m_n}^n:= y_n$ and define $x_{m_n-1}^n:= g_{m_n}hg_{m_n}^{-1}g_{m_n-1}$. Given the point $x_{m_n-j}^n$, where $j\in \{1, 2, \ldots, m_n - 1\}$, set $x_{m_n-j-1}^n:= x_{m_n-j}^ng_{m_n-j}^{-1}g_{m_n-j-1}$. Note that for each $i$, $x_i^n\in X_i = g_i\Gamma_H$, and so $\underline{x}^n$ defined in this fashion is a 1-chain in $X(\gamma)$ with terminal point $y_n$.

We now have a sequence of 1-chains $\underline{x}^n$ with terminal points $y_n$ converging to $y\in \partial X(\gamma)$. Passing to a subsequence, we may assume that $x_0^n$ converges to $x_0\in X_0\cup \partial X_0$. Suppose first that $x_0\in X_0$. Then, since the points $x_1^n$ remain in a compact subset of $X_1$, they subconverge on a point $x_1\in X_1$ with $d_{X(\gamma)}(x_0, x_1) = 1$. Continuing on in this fashion, we can pass to a subsequence of our partial chains to get an infinite 1-chain $\underline{x}=\{x_0, x_1, \ldots, x_n, \ldots\}$ in $X(\gamma)$, where $x_i^n$ converges to $x_i\in X_i$ for all $i$. Note that for large enough $n$, $x_i^n$ remains uniformly close to $x_i$ for arbitrarily many $i$. Hence, we must have that the terminal points of the chains $\underline{x}^n$ converge to the terminal point of $\underline{x}$ in $X(\gamma)\cup \partial X(\gamma)$. Since the chains $\underline{x}^n$ each have terminal point $y_n$, we therefore have that $y\in \partial X(\gamma)$ is the terminal point of a 1-chain in $X(\gamma)$, and so $y\in \partial^+ X(\gamma)$. 
Suppose now that $x_0\in \partial X_0$. Then, by \Cref{Bowditch Lemma 2.3.5} (1), we have that $\partial \omega_\gamma(x_0) = y$, where $\partial \omega_\gamma: \partial X_0\to \partial X(\gamma)$. Therefore, for all $y\in \partial X(\gamma)$, either $y$ is the endpoint of a 1-chain in $X(\gamma)$, or $y\in \omega_\gamma(\partial X_0)$. 

So, suppose that $(x_i)_{i\in\mb{Z}}$ is a 1-chain. We have that for all $i, j\in \mb{Z}$ with $i < j$, \begin{align*}
d_Q(Px_i, Px_j) &\leq d_{X(\gamma)}(x_i, x_j) \\
&\leq d_{X(\gamma)}(x_i, x_{i+1}) + d_{X(\gamma)}(x_{i+1}, x_{i+1}) + \cdots + d_{X(\gamma)}(x_{j-1}, x_{j}) \\
&= d_Q(Px_i, Px_{i+1}) + \cdots + d_Q(Px_{j-1}, Px_j) \\
&= d_Q(z_i, z_j).
\end{align*} Hence, every 1-chain in $X(\gamma)$ interpolates a geodesic in $X(\gamma)$ which is an isometric lift of $\gamma$. Furthermore, as for all $i,j\in \mb{Z}$, $d_Q(Px_i, Px_j)\leq d_G(x_i, x_j)\leq d_{X(\gamma)}(x_i, x_j)$, we have that every 1-chain interpolates a geodesic in $\Gamma_G$ as well. Therefore, if $y\in \partial^+ X(\gamma)$ is the terminal point in $X(\gamma)$ of the positive 1-chain $(y_i)_{i\in\mb{N}}$, then the terminal point of this chain in $\Gamma_G$ will determine a point of $\partial G$ as well. 
As the only $r$-chains we will be considering in $X(\gamma)$ are 1-chains, all 1-chains in this space will now simply be referred to as chains.

\begin{convention}\label{StackInclusionConvention}
Given $P: \Gamma_G\to \Gamma_Q$ as in \Cref{ShortExactSequenceConvention} and $\gamma$ as in \Cref{StackConvention}, let $\sigma: \gamma\to \Gamma_G$ denote an isometric lift of $\gamma$ such that for all $z_i\in \gamma$, $P(\sigma(z_i)) = z_i$ and set $g_i:= \sigma(z_i)$. Let $X(\gamma):= P^{-1}(\gamma)$ and $X(\gamma)^+:= P^{-1}(\gamma^+)$ be the stacks which consist of the sheets $X_i = g_i \Gamma_H$ for all $i\in \mb{Z}$ and $i\in \mb{N}$, respectively. 
Denote by $\omega_\gamma: X_0\to X(\gamma)$ and $\omega_\gamma^+: X_0\to X(\gamma)^+$ the inclusions of the sheet $X_0 = g_0\Gamma_H$ into $X(\gamma)$ and $X(\gamma)^+$ respectively. 
Define $i_{X_0}: \Gamma_H\to X_0$ as follows. Set $i_{X_0}(h) := g_0 \cdot g_0^{-1} h g_0 = h g_0$ for all vertices $h\in \Gamma_H$. Extend $i_{X_0}$ to a map on all of $\Gamma_H$ by sending an edge $[a,b]$ to a shortest path between $ag_0$ and $bg_0$. 
Now let $i_\gamma: \Gamma_H\to X(\gamma)$ be given by $i_\gamma:= \omega_\gamma \circ i_{X_0}$ and $i_\gamma^+: \Gamma_H \to X(\gamma)^+$ be given by $i_\gamma^+:= \omega_\gamma^+\circ i_{X_0}$. Note that if $1\in \gamma$, then $i_{X_0}$ and $i_\gamma^+$ are simply the identity inclusion map $i: \Gamma_H\to \Gamma_G$.
\end{convention}

\begin{lemma}\label{i_gamma Cannon-Thurston}
The maps $i_\gamma: \Gamma_H\to X(\gamma)$ and $i_\gamma^+: \Gamma_H\to X(\gamma)^+$ as given in \Cref{StackInclusionConvention} extend continuously to the maps $\widehat{i}_\gamma: \widehat{\Gamma_H} \to \widehat{X(\gamma)}$ and $\widehat{i}_\gamma^+: \widehat{\Gamma_H} \to \widehat{X(\gamma)^+}$, respectively.
\end{lemma}
\begin{proof}
Given $\omega_\gamma: X_0\to X(\gamma)$ and $\omega_\gamma^+: X_0\to X(\gamma)^+$ as in \Cref{StackInclusionConvention}, note that \Cref{Bowditch Prop 2.3.2} gives that the Cannon-Thurston maps $\partial \omega_\gamma: \partial X_0\to \partial X(\gamma)$ and $\partial \omega_\gamma^+: \partial X_0 \to \partial X(\gamma)$ both exist. Let $\widehat{\omega_\gamma}: \widehat{X_0}\to \widehat{X(\gamma)}$ and $\widehat{\omega_\gamma^+}: \widehat{X_0}\to \widehat{X(\gamma)^+}$ denote the continuous extensions of $\omega_\gamma$ and $\omega_\gamma^+$. For all $g\in G$, conjugation by $g$ gives an automorphism of $H$ which takes $h\in H$ to $g^{-1}hg$. This automorphism is a quasi-isometry from $\Gamma_H$ to itself. So, $i_{X_0}: \Gamma_H\to X_0$ is a quasi-isometry from $\Gamma_H$ to $X_0 = g_0 \Gamma_H$, and so extends to a homeomorphism $\partial i_{X_0}: \partial H \to \partial X_0$. Hence, $\partial i_\gamma:= \partial \omega_\gamma \circ \partial i_{X_0}$ and $\partial i_\gamma^+: \partial \omega_\gamma^+\circ \partial i_{X_0}$ exist, are continuous, and extend $i_\gamma$ and $i_\gamma^+$ continuously to the maps $\widehat{i_\gamma}$ and $\widehat{i_\gamma^+}$, respectively.
\end{proof}

\Cref{i_gamma Cannon-Thurston} allows us to now refer to the maps $\partial i_\gamma$ and $\partial i_\gamma^+$ as Cannon-Thurston maps. The goal of the remainder of this section is to show that the maps $\partial i_\gamma$ and $\partial i_\gamma^+$ are surjective. We will first show surjectivity for the case where the geodesic $\gamma$ lives over the identity in $\Gamma_Q$.

\begin{convention}\label{IdentityStackConvention}
Let $\gamma = (z',z)$ be as in \Cref{StackConvention} and let $\gamma': = z_0^{-1}\cdot \gamma = (z_0^{-1}z', z_0^{-1}z)$. Note that $1\in V(\gamma')$. For each $z_i\in \gamma$, let $z_i':= z_0^{-1}\cdot z_i$. 
Given $\sigma: \gamma\to \Gamma_G$ as in \Cref{StackInclusionConvention},
let $\sigma': \gamma'\to \Gamma_G$ be such that $\sigma':= g_0^{-1}\cdot \sigma$. Set $g_i':= \sigma'(z_i')$, and denote the sheet $g_i'\Gamma_H$ by $X_i'$. Note that the sheet $X_0'$ is the identity coset $1\cdot\Gamma_H$, and so the map $i_{X_0'}: \Gamma_H\to X_0'$ is the identity map.
\end{convention}

In a similar manner as Bowditch \cite{Bow13}, we define a map $\widehat{\tau}_{\gamma'}: \widehat{X(\gamma')} \rightarrow \widehat{\Gamma_G}$ with $\widehat{i} = \widehat{\tau}_{\gamma'} \circ \widehat{i}_{\gamma'}$ and will later show that $\partial \tau_{\gamma'}: \partial X(\gamma')\to \partial G$ is continuous. 
Let $\tau_{\gamma'} := \widehat{\tau}_{\gamma'} |_{X(\gamma')}$ be the identity inclusion of $X(\gamma')$ into $\Gamma_G$ given by $\tau_{\gamma'}(g) =g$. 
Note that for all $h\in H$, $\tau_{\gamma'} \circ i_{\gamma'} (h) = h = i(h)$. 

As the map $i_{X_0'}$ is the identity map, $\partial \omega_{\gamma'} = \partial i_{\gamma'}$. So by \Cref{Bowditch Prop 2.3.2}, we have that $\partial X(\gamma') = \partial i_{\gamma'}(\partial H) \cup \partial^\pm X(\gamma')$. 
If $(y_i)_{i\in \mb{N}}$ is a positive $1$-chain in $X(\gamma')$ with endpoint $y\in \partial^+X(\gamma')$, then $(\tau_{\gamma'}(y_i))_{i\in\mb{N}}$ interpolates a geodesic ray in $\Gamma_G$ with the same label as the geodesic ray interpolated by $(y_i)$ in $X(\gamma')$. Denote the endpoint of this geodesic ray in $\Gamma_G$ by $\overline{y}\in \partial G$, and for all $y\in \partial^\pm X(\gamma')$ define $\partial \tau_{\gamma'}(y) := \overline{y}$. Finally, for all $a\in \partial H$, define $\partial \tau_{\gamma'}(\partial i_{\gamma'}(a)) := \partial i(a)$. 
Note that if $(x_i)_{i\in \mathbb{N}}$ and $(y_i)_{i\in\mathbb{N}}$ are distinct but equivalent $1$-chains in $X(\gamma')$, then the geodesic rays interpolated by these chains are Hausdorff close in both $X(\gamma')$ and $\Gamma_G$. Hence, $\tau_{\gamma'}$ is well-defined on equivalence classes of chains. 
To finish showing that $\widehat{\tau}_{\gamma'}$ is well-defined, we need the following lemma.

\begin{lemma}\label{Analogue of Bowditch Lemma 2.3.1}
Let $\gamma'$ be as in \cref{IdentityStackConvention}. Suppose $(\underline{x}^n)_{n\in \mb{N}}$ is a sequence of positive chains in $X(\gamma')$, where $\underline{x}^n = (x_i^n)_{i\in \mb{N}}$ is a positive chain with terminal point $y_n\in \partial^+ X(\gamma')$. Suppose also that in $\widehat{X(\gamma')}$, $y_n\rightarrow y\in \partial X(\gamma')$ and in $\widehat{X_0'}$, $x_0^n\rightarrow \partial i_{X_0'}(a)\in \partial X_0'$. Then in $\widehat{\Gamma_G}$, $\widehat{\tau}_{\gamma'}(y_n)\rightarrow \widehat{i}(a)$.
\end{lemma}

\begin{proof}
Let $f(n) = \max\{d_{X_0'}(a,b) \mid d_G(\tau_{\gamma'}(a),\tau_{\gamma'}(b))\leq n\}$. Note that since $\Gamma_G$ is finitely generated, such a maximum exists and that $f(n)\rightarrow \infty$ as $n\rightarrow \infty$. For each $n\in \mb{N}$, there exists $a_n\in \Gamma_H$ such that $x_0^n = i_{X_0'}(a_n) = a_n$. As $x_0^n\to \partial i_{X_0'}(a)$, this implies that $a_n\to a\in \partial H$ in $\widehat{\Gamma_H}$. Let $\lambda_n = [\widehat{\tau}_{\gamma'}(x_0^n), \widehat{\tau}_{\gamma'}(y_n))_G = [i(a_n), \widehat{\tau}_{\gamma'}(y_n))_G$ be the geodesic ray in $\Gamma_G$ interpolated by $(\tau_{\gamma'}(x_i^n))_{i\in\mb{N}}$ for each $n\in \mb{N}$.

Suppose that in $\widehat{\Gamma_G}$, $\lim_{n\to \infty} \widehat{\tau}_{\gamma'}(y_n) \neq \lim_{n\to\infty} i(a_n)$. Then, there exist constants $R,N > 0$ such that for all $n\geq N$, $d_G(1,\lambda_n)\leq R$. So, for each $n\geq N$ there exists some point $x_{i_n}^n$ in the chain $\underline{x}^n$ such that $d_G(1,\tau_{\gamma'}(x_{i_n}^n)) \leq R$. Then, we have that \begin{align*}
    d_G(1, \tau_{\gamma'}(x_{i_n}^n)) & \geq d_Q(P\cdot 1, P\cdot \tau_{\gamma'}(x_{i_n}^n)) \\
    &= d_Q(1, P\cdot x_{i_n}^n) \\
    &= |i_n|.
\end{align*} As $d_G(1,\tau_{\gamma'}(x_{i_n}^n))\leq R$, this means that $|i_n|\leq R$. Note that since $\underline{x}^n$ is a 1-chain, we have that $d_G(\tau_{\gamma'}(x_0^n), \tau_{\gamma'}(x_{i_n}^n)) = d_{X(\gamma')}(x_0^n, x_{i_n}^n) = |i_n|\leq R$ So,
\begin{align*}
d_{X_0'}(1, x_0^n)&\leq f(d_G(\tau_{\gamma'}(1), \tau_{\gamma'}(x_0^n))) \\
&= f(d_G(1, \tau_{\gamma'}(x_0^n))) \\
&\leq f(d_G(1, \tau_{\gamma'}(x_{i_n}^n)) + d_G(\tau_{\gamma'}(x_{i_n}^n), \tau_{\gamma'}(x_0^n))) \\
&\leq f(2R).
\end{align*} 

But, $d_{X_0'}(1, x_0^n)\rightarrow \infty$ as $n\rightarrow \infty$ since $x_0^n\to \partial i_{X_0'}(a)\in \partial X_0'$, and so we have a contradiction. Therefore, $d_G(1, \lambda_n)\rightarrow\infty$ as $n\rightarrow\infty$. Hence, in $\widehat{\Gamma_G}$, $\lim_{n\rightarrow \infty} \tau_{\gamma'}(x_0^n) = \lim_{n\rightarrow \infty} \widehat{\tau}_{\gamma'}(y_n)$. As $\tau_{\gamma'}(x_0^n) = i(a_n)$ and $i(a_n)\to \widehat{i}(a)$ as $n\to \infty$, we have that $\widehat{\tau}_{\gamma'}(y_n)\rightarrow \widehat{i}(a)$ as desired.
\end{proof}

\begin{lemma}\label{TauWellDefined}
The map $\widehat{\tau}_{\gamma'}: \widehat{X(\gamma')}\to \widehat{\Gamma_G}$ is well-defined and satisfies $\widehat{\tau}_{\gamma'}\circ \widehat{i}_{\gamma'} = \widehat{i}: \widehat{\Gamma_G}\to \widehat{\Gamma_H}$
\end{lemma}

\begin{proof}
If $x\in \Gamma_H$, then $\tau_{\gamma'} \circ i_{\gamma'}(x) = x = i(x)$. Similarly, if $a\in \partial H$, then $\partial \tau_{\gamma'} \circ \partial i_{\gamma'}(a) = \partial i(a)$. So, $\widehat{i} = \widehat{\tau}_{\gamma'} \circ \widehat{i}_{\gamma'}$. Now, it suffices to show that $\partial \tau_{\gamma'}: \partial X(\gamma')\to \partial G$ is well-defined. 
First, we need to show that if $y\in \partial^+ X(\gamma')$ and $a\in \partial H$ are such that $\partial i_{\gamma'}(a) = y$, then $\partial \tau_{\gamma'}(y) = \partial \tau_{\gamma'}(\partial i_{\gamma'}(a))$. So, suppose that $y\in \partial^+ X(\gamma')$ and $a\in \partial H$ are such that $\partial i_{\gamma'}(a) = y$. Since $\partial i_{\gamma'}(a) = y$ and $\partial i_{X_0'}$ is the identity,
this implies that $\partial i_{\gamma'}(a) = \partial \omega_{\gamma'} \circ \partial i_{X_0'}(a) = \partial \omega_{\gamma'}(a) = y$.  By \Cref{Bowditch Lemma 2.3.5} (1), there exists a sequence $(\underline{x}^n)_n$ of positive chains, each converging to $y$, with $x_0^n$ converging to $a = \partial i_{X_0'}(a)\in \partial X_0'$. By \Cref{Analogue of Bowditch Lemma 2.3.1}, the existence of such a sequence of chains implies that $\partial \tau_{\gamma'}(y) = \partial i(a)$. Hence, $\partial \tau_{\gamma'}(y) = \partial \tau_{\gamma'}(\partial i_{\gamma'}(a))$.

Now, suppose that $a,b\in \partial H$ with $a\neq b$ are such that $\partial i_{\gamma'}(a) = \partial i_{\gamma'}(b)$. Since $\partial i_{X_0'}$ is the identity, this implies that $\partial \omega_{\gamma'}(a) = \partial \omega_{\gamma'}(b)$. By \Cref{Bowditch Prop 2.3.9} (4), we may assume without loss of generality that $\partial \omega_{\gamma'}^+(a) = \partial \omega_{\gamma'}^+(b)$. Since $a$ and $b$ are distinct, we have by \Cref{Bowditch Prop 2.3.7} (3) that $\partial \omega_{\gamma'}^+(a) = \partial \omega_{\gamma'}^+(b) = y\in \partial^+X(\gamma')$. By \Cref{Bowditch Prop 2.3.6} (2), we now have that $\partial \omega_{\gamma'}(a) = \partial \omega_{\gamma'}(b) = y$. So by the same reasoning as above, \Cref{Bowditch Lemma 2.3.5} (1) and \Cref{Analogue of Bowditch Lemma 2.3.1} give that $\partial \tau_{\gamma'}(\partial i_{\gamma'}(a)) = \partial \tau_{\gamma'}(\partial i_{\gamma'}(b)) = \partial \tau_{\gamma'}(y) = \overline{y}$.
\end{proof}

\begin{corollary}\label{igamma=i}
If $a,b\in \partial H$ are such that $\partial i_{\gamma'}^+(a) = \partial i_{\gamma'}^+(b)$ then $\partial i(a) = \partial i(b)$.
\end{corollary}

\begin{proof}
Suppose $a,b\in \partial H$ are such that $\partial i_{\gamma'}^+(a) = \partial i_{\gamma'}^+(b)$. If $a = b$, then $\partial i(a) = \partial i(b)$. So, suppose $a\neq b$. 
As $\partial i_{X_0'}$ is the identity, we have that $\partial \omega_{\gamma'}^+(a) = \partial \omega_{\gamma'}^+(b)$. By \Cref{Bowditch Prop 2.3.7} (3), there exists $y\in \partial^+ X(\gamma')$ such that $\partial \omega_{\gamma'}^+(a) = \partial \omega_{\gamma'}^+(b) = y$. By \Cref{Bowditch Prop 2.3.7} (2), this implies that $\partial \omega_{\gamma'}(a) = \partial \omega_{\gamma'}(b)$. So, $\partial i_{\gamma'}(a) = \partial \omega_{\gamma'} \circ \partial i_{X_0'}(a) = \partial \omega_{\gamma'} \circ \partial i_{X_0'}(b) = \partial i_{\gamma'}(b)$. As $\partial \tau_{\gamma'}$ is well-defined by \Cref{TauWellDefined}, we have that $\partial \tau_{\gamma'}(\partial i_{\gamma'}(a)) = \partial \tau_{\gamma'}(\partial i_{\gamma'}(b))$, and so $\partial i(a) = \partial i(b)$. 
\end{proof}

The goal of the remainder of this section is to use this work of Bowditch to prove that the Cannon-Thurston map $\partial i_\gamma^+: \partial H \rightarrow \partial X(\gamma)^+$ is surjective.

\begin{lemma}\label{QuasigeodesicLemma}
Fix $\gamma = (z',z)\subseteq \Gamma_Q$ as in $\Cref{StackConvention}$ and let $X(\gamma)^+$ be as described above. Let $(y_n)_{n\in\mb{N}}$ be a $1$-chain in $X(\gamma)^+$ and denote the word which labels the geodesic from $y_0$ to $y_n$ in $X(\gamma)^+$ by $\alpha_n$. Fix some $h\in H$ of infinite order and let $\rho_n$ denote any path in $X(\gamma)^+$ which is the concatenation of a path labeled by $\alpha_n$ followed by a path labeled by $h$ and finally a path labeled by $\alpha_n^{-1}$. Then, there exists some constant $C\geq 0$ independent of $n$ (but dependent on $h$) such that for all $n$, $\rho_n$ is a $(1,C)$-quasigeodesic in $X(\gamma)^+$.
\end{lemma}

\begin{proof}
Let $(y_n)_{n\in\mb{N}}$ be a 1-chain in $X(\gamma)^+$, and for each $n\geq 0$ let $\alpha_n$ denote the word which labels the geodesic from $y_0$ to $y_n$. Given $h\in H$, let $\beta$ denote any quasigeodesic in $X(\gamma)^+$ labeled by $h$. Since $h\in H$ is fixed, there exists some constant $C'\geq 0$  such that $\beta$ is a $(1, C')$-quasigeodesic. For each $n\geq 0$, let $[x_n = y_0, y_n]$ be the geodesic in $X(\gamma)^+$ from $x_n = y_0$ to $y_n$ labeled by $\alpha_n$, let $z_n\in X(\gamma)^+$ be a point such that $\beta$ is a quasigeodesic in $X(\gamma)^+$ from $y_n$ to $z_n$, and let $[z_n, w_n]$ be the geodesic in $X(\gamma)^+$ labeled by $\alpha_n^{-1}$. Denote by $\delta_n$ the label of the geodesic in $X(\gamma)^+$ between $x_n$ and $w_n$.

For each $n\geq 0$, consider the quadrilateral in $X(\gamma)^+$ with vertices $y_0 = x_n$, $y_n$, $z_n$, $w_n$, and with sides $[x_n, y_n]$ labeled by $\alpha_n$, $\beta$ labeled by $h$, $[z_n, w_n]$ labeled by $\alpha_n^{-1}$, and $[x_n, w_n]$ labeled by $\delta_n$.
Unless otherwise specified, we will denote $d_{X(\gamma)^+}$ simply by $d$, and all geodesic and quasigeodesic segments considered are geodesics or quasigeodesics in $X(\gamma)^+$.  

As before, we need to show that if $p$ and $q$ are arbitrary points on $\rho_n = [x_n, y_n]\cup \beta \cup [z_n, w_n]$, then the distance between $p$ and $q$ along $\rho_n$ is at most $d(p,q) + C$. There are two cases to consider. By \Cref{ThinQuadrilateral}, either there is a point on $[x_n,w_n]$ at most distance $2\delta$ in $X(\gamma)^+$ from a point on $[y_n,z_n]$, or there is a point on the side $[x_n,y_n]$ at most distance $2\delta$ in $X(\gamma)^+$ from a point on the side $[z_n,w_n]$. If there is some point on the side $[x_n,w_n]$ within $2\delta$ of a point on the side $[y_n,z_n]$, then \Cref{GeodesicQuadrilateral} gives that $[x_n, y_n]\cup [y_n, z_n] \cup [z_n, w_n]$ is a $(1, 4\delta + 4d(y_n, z_n))$-quasigeodesic. Since $\beta$ is a $(1,C')$-quasigeodesic between $y_n$ and $z_n$, this gives that $\rho_n$ is a $(1, C)$-quasigeodesic for some $C\geq 0$.

So, suppose now that the two sides labeled by $\alpha_n$ and $\alpha_n^{-1}$ come within $2\delta$ of each other in $X(\gamma)^+$. We make the following claim: 

\textbf{Claim:} If $a\in [x_n,y_n]$ and $a'\in [z_n,w_n]$ are the furthest points in $X(\gamma)^+$ from $y_n$ and $z_n$, respectively, such that $d(a,a')\leq 2\delta$, then there is some constant $K > 0$ dependent on $h$ but independent of $n$ such that $\max\{d(a,y_n), d(a',z_n)\}\leq K$.

Assuming this claim, we will now show that $\rho_n$ is a $(1,C)$-quasigeodesic in $X(\gamma)^+$. First fix $p\in [x_n,y_n]$ and $q\in \beta$. Since in $X(\gamma)^+$, $(p,q; X(\gamma)^+)_{y_n}$ is bounded by $|\beta|_{X(\gamma)^+}\leq |h|_H$, we have that \begin{align*}
    d(p,y_n) + d(y_n,q) &= d(p,q) + 2(p,q; X(\gamma)^+)_{y_n} \\&\leq d(p,q) + 2|h|_H.
\end{align*}

So, suppose $p\in [x_n,y_n]$ and $q\in [z_n,w_n]$. If $p\in [a,y_n]$ and $q\in [a',z_n]$, then $d(p,y_n) + |\beta|_{X(\gamma)^+}+ d(z_n,q) \leq d(p,q) + |h|_{H} + 2K$. Now suppose $p\in [x_n,a]$ and $q\in [a',z_n]$. Since $d(a,a')\leq 4\delta$, we have by the triangle inequality that \begin{align*}
    d(q,a) &\leq d(q,a') + d(a',a) \\ &\leq K + 2\delta, \text{ and}
\end{align*}
\begin{align*}
    d(p,a) &\leq d(p,q) + d(q,a).
\end{align*}
Therefore, \begin{align*}
    d(p,a) + d(a,y_n) + |\beta|_{X(\gamma)^+} + d(z_n,q) &\leq d(p,q) + d(q,a) + K + |h|_{H} + K \\ &\leq d(p,q) + 3K + 2\delta + |h|_{H}.
\end{align*}
The final case to consider is when $p\in [x_n,a]$ and $q\in [w_n,a']$. In this case, there must be a point $u\in [p,q]$ and $v\in [a,a']$ such that $d(u,v)\leq 2\delta$. This is because by choice of $a$ and $a'$, there are no points at which $[q,a']$ is within a distance of $2\delta$ of $[p,a]$ in $X(\gamma)^+$. So, we have that $d(q,v) \leq d(q,u) + d(u,v)$ and $d(p,v)\leq d(p,u) + d(u,v)$. Additionally, $d(p,a) \leq d(p,v) + d(v,a)$ and $d(q,a')\leq d(q,v) + d(v,a')$. Hence, we have that \begin{align*}
    d(p,a) + d(a,y_n) + |\beta|_{X(\gamma)^+} + d(z_n,a') + d(a',q) &\leq d(p,a) + K + |h|_H + K + d(a',q) \\ 
    &\leq d(p,v) + d(v,a) + 2K + |h|_H +  d(q,v) + d(v,a') \\
    &\leq d(p,v) + d(q,v) + 2K + |h|_H + 2\delta \\
    &\leq d(p,u) + d(q,u) + 2d(u,v) + 2K + |h|_H + 2\delta \\
    &\leq d(p,q) + 2K + |h|_H + 6\delta.
\end{align*}

Proof of Claim: Suppose to the contrary that there is no such bound on the how long the sides labeled by $\alpha_n$ and $\alpha_n^{-1}$ stay uniformly close in $X(\gamma)^+$. Let $S_Q$ be the generating set for $Q$ and let $L = \{w\in \Sigma_Q^* \mid w \text{ a geodesic in }Q\}$. Since $Q$ is a hyperbolic group, the language $L$ of geodesic words is a regular language for $Q$ (see \cite{EpsEtAl92}) which is accepted by some finite state automaton, $\mathcal{A}$, with start state $s_0$. Then, $\gamma^+ = [z_0,z)\subseteq \Gamma_Q$ gives an infinite path from $s_0$ in $\mathcal{A}$ such that all states are accept states. Let $\gamma_n$ denote the initial portion of the path $\gamma^+$ of length $n$, i.e., $\gamma_n:= P([y_0 = x_n,y_n])$. 

For each $n$, assume without loss of generality that the side of $\rho_n$ labeled by $\alpha_n$ begins at the vertex $y_0$ and ends at the vertex $y_n$. Let $y_{i_n}$ denote the vertex along the side $\alpha_n$ where the side labeled by $\alpha_n$ and the side labeled by $\alpha_n^{-1}$ begin to be $2\delta$ close. Note that after the point 
$y_{i_n}$, the sides labeled by $\alpha_n$ and $\alpha_n^{-1}$ will continue to travel within a distance of $|h|_{X(\gamma)^+}$ of each other in $X(\gamma)^+$. Project the $X(\gamma)^+$-geodesic $[y_{i_n}, y_n]$ to $Q$ and feed this geodesic, $P([y_{i_n}, y_n])$, into $\mathcal{A}$. Note that by assumption, the length of these geodesics go to infinity as $n\to \infty$. So, there will be some $n>0$ for which some state in $\mathcal{A}$ repeats more times than the number of words in $G$ of length at most $|h|_{X(\gamma)^+}$. Note that the label of any loop in $\mathcal{A}$ is a periodic $Q$-geodesic word. 
Since there is a state that repeats more times than the number of words in $G$ of length at most $|h|_{X(\gamma)^+}$, it follows that there is some subpath of $[y_{i_n}, y_n]$ labeled by a word $v\in \Sigma^*_Q$ which has infinite order in $Q$ and some word $m\in \Sigma_G^*$ of length at most $|h|_{X(\gamma)^+}$ such that in $G$, $P^{-1}(v)m(P^{-1}(v))^{-1} = m$ and such that $h$ is conjugate to $m$ in $G$. As $h$ has infinite order in $G$ and $h$ is conjugate to $m$, it follows that $m$ is infinite order in $G$ as well. As $P^{-1}(v)$ and $m$ commute in $G$, this implies that $(P^{-1}(v))^p = m^q$, for some $p,q\neq 0$. But then $v^p = 1$ in $Q$, because $h$ projects to the identity in $Q$ which means that $m$ projects to the identity in $Q$ as well. The fact that $v^p = 1$ contradicts $v$ being a periodic geodesic in $Q$. This completes the proof of the claim and the lemma.
\end{proof}

\begin{mainthrmB}\hypertarget{thrmB}{}
Let $1\to H\to G\to Q\to 1$ be a short exact sequence of infinite, finitely generated, word-hyperbolic groups. Let $z,z'\in \partial Q$ be distinct and let $\gamma\subseteq \Gamma_Q$ be a bi-infinite geodesic in $\Gamma_Q$ between $z$ and $z'$. Let $i_\gamma^+: \Gamma_H\to X(\gamma)^+$ be the inclusion of $\Gamma_H$ into the semi-infinite stack $X(\gamma)^+$ over $\gamma^+ = [z_0,z)$, and let $i_\gamma: \Gamma_H\to X(\gamma)$ be the inclusion of $\Gamma_H$ into the bi-infinite stack $X(\gamma)$, as in \Cref{StackInclusionConvention}. Then,
\begin{enumerate}
    \item the Cannon-Thurston map $\partial i_\gamma^+: \partial H\to \partial X(\gamma)^+$ is surjective; and
    \item the Cannon-Thurston map $\partial i_\gamma: \partial H\to \partial X(\gamma)$ is surjective.
\end{enumerate}
\end{mainthrmB}

\begin{proof}
Let $\gamma = (z',z)\subseteq \Gamma_Q$ be as in \cref{StackConvention} and let $\gamma':= z_0^{-1}\cdot \gamma$ be as in \cref{IdentityStackConvention}. We will first show that the Cannon-Thurston maps $\partial i_{\gamma'}^+: \partial H\to \partial X(\gamma')^+$ and $\partial i_{\gamma'}: \partial H\to \partial X(\gamma')$ are surjective.

Consider first the map $\partial i_{\gamma'}^+: \partial H\to \partial X(\gamma)^+$. Since $\partial i_{\gamma'}^+ = \partial\omega_{\gamma'}^+ \circ \partial i_{X_0'}$ and $\partial i_{X_0'}$ is the identity, it suffices to show that $\partial \omega_{\gamma'}^+$ is surjective. By \Cref{Bowditch Prop 2.3.2} (3), we need only show that if $y\in \partial^+X(\gamma')^+$, then there exists $a\in \partial X_0'$ such that $\partial i_{\gamma'}^+(a) = y$. So, suppose that $y\in \partial^+ X(\gamma')^+$ is the endpoint of the chain $(y_n)$ and fix some $h\in H$ of infinite order. Let $\alpha_n$ be the word which labels the path from $y_0$ to $y_n$ in $X(\gamma')^+$, and consider the path $\rho_n$ in $X(\gamma')^+$ which is labeled by the word $\alpha_n h \alpha_n^{-1}$. 

By \Cref{QuasigeodesicLemma}, $\rho_n$ is a $(1, C)$-quasigeodesic in $X(\gamma')^+$ for some $C$ independent of $n$. Let $h_n$ be the word which labels the geodesic in $X_0'$ between the endpoints of $\rho_n$. Since $|h_n|_H \to \infty$, there exists a subsequence $h_{n_i}$ such that $y_0h_{n_i}\to a\in \partial X_0'$. Since $\partial \omega_{\gamma'}^+$ is a continuous extension of $\omega_{\gamma'}^+$, we have that 
$$\lim_{n_i\to \infty} \omega_{\gamma'}^+(y_0h_{n_i}) = \partial \omega_{\gamma'}^+ \lim_{n_i\to \infty} y_0h_{n_i} = \partial \omega_{\gamma'}^+(a).$$
Since $y_{n_i}\to y$ and since $\rho_{n_i}$ is a quasigeodesic and $y_{n_i}\in \rho_{n_i}$, it follows that $\lim_{n_i\to \infty} y_{n_i} = \lim_{n_i\to \infty} \omega_\gamma^+(y_0h_{n_i}) = y$ in $\widehat{X(\gamma')^+}$.

To see that $\partial i_{\gamma'}: \partial H\to \partial X(\gamma')$ is surjective, note that by \Cref{Bowditch Prop 2.3.3} (2), $\partial X(\gamma') = \partial^+X(\gamma') \cup \partial^- X(\gamma') \cup \partial i_{\gamma'}(\partial H)$. Note that the map $\widehat{i}_{\gamma'}: \widehat{H}\to \widehat{X(\gamma')}$ is defined in the same way as $\widehat{i}_{\gamma'}^+$. So, to show the surjectivity of $\partial i_{\gamma'}$, it suffices to note that in the above argument, we can replace $y\in \partial^+X(\gamma')^+$ with $y'\in \partial^-X(\gamma')^-$. As the same reasoning holds, we have that $\partial i_{\gamma'}: \partial H\to \partial X(\gamma')$ is surjective as well. 

Now, let $t_{g_0}^H: \Gamma_H\to g_0\Gamma_H$, $t_{g_0}: X(\gamma')\to X(\gamma)$, and $t_{g_0}^+: X(\gamma')^+\to X(\gamma)^+$ denote the maps induced by left-translation of the vertices of $\Gamma_H$, $X(\gamma')$, and $X(\gamma')^+$, respectively, by the element $g_0 = \sigma(z_0)$. 
Note that for all $h\in H$, $\omega_\gamma \circ t_{g_0}^H(h) = t_{g_0}\circ i_{\gamma'}(h)$ and $\omega_{\gamma}^+ \circ t_{g_0}^H(h) = t_{g_0}^+\circ i_{\gamma'}^+(h)$. Since $t_{g_0}^H$, $t_{g_0}$, and $t_{g_0}^+$ are isometries, these maps extend continuously to the boundary maps $\partial t_{g_0}^H: \partial H\to \partial g_0H$, $\partial t_{g_0}: \partial X(\gamma')\to \partial X(\gamma)$, and $\partial t_{g_0}^+: \partial X(\gamma')^+\to \partial X(\gamma)$, respectively, which are homeomorphisms. Hence, we have that for all $a\in \partial H$, $\partial \omega_\gamma\circ \partial t_{g_0}^H(a) = \partial t_{g_0}\circ \partial i_{\gamma'}(a)$ and $\partial \omega_\gamma^+\circ \partial t_{g_0}^H(a) = \partial t_{g_0}^+\circ \partial i_{\gamma'}^+(a)$. As $\partial i_{\gamma'}$ and $\partial i_{\gamma'}^+$ are surjective by the above argument and as $\partial t_{g_0}^H$, $\partial t_{g_0}$, and $\partial t_{g_0}^+$ are homeomorphisms, this implies that $\partial \omega_\gamma$ and $\partial \omega_\gamma^+$ are surjective.

As noted previously, each $g\in G$ gives rise to an automorphism $\phi_g$ of $H$ with $\phi_g(h) = g^{-1}hg$. This automorphism of $H$ induces a quasi-isometry of $\Gamma_H$ taking an edge $[u,v]$ to a shortest edge path between $\phi_g(u)$ and $\phi_g(v)$. 
As $\phi_{g}: \Gamma_H\to \Gamma_H$ is a quasi-isometry, it extends to a homeomorphism $\partial \phi_{g}: \partial H\to \partial H$. Recall that $i_\gamma = \omega_\gamma \circ t_{g_0}^H\circ \phi_{g_0}$ and $i_\gamma^+ = \omega_\gamma^+ \circ t_{g_0}^H\circ \phi_{g_0}$. So, $\partial i_\gamma = \partial \omega_\gamma \circ \partial t_{g_0}^H \circ \partial \phi_{g_0}$ and $\partial i_\gamma^+ = \partial \omega_\gamma^+ \circ \partial t_{g_0}^H \circ \partial \phi_{g_0}$. As $\partial \omega_\gamma$ and $\partial \omega_\gamma^+$ are surjective, and as $\partial t_{g_0}^H$ and $\partial \phi_{g_0}$ are homeomorphisms, we have that $\partial i_\gamma$ and $\partial i_\gamma^+$ are surjective.
\end{proof}

Recall that given the maps $\partial i_\gamma^+: \partial H\to \partial X(\gamma)^+$ and $\partial i_\gamma: \partial H\to \partial X(\gamma)$, Bowditch defines a map $\partial \tau^+: \partial X(\gamma)^+\to \partial X(\gamma)$ with $\partial i_\gamma = \partial \tau^+ \circ \partial i_\gamma^+$. This map is given by $\partial \tau^+(y) = y$ for all $y\in \partial^+ X(\gamma)^+$, and $\partial \tau^+ \circ \partial i_\gamma^+(a) = \partial i_\gamma(a)$ for all $a\in \partial H$. We can now show the following about the map $\partial \tau^+$.

\begin{corollary}
The map $\partial \tau^+: \partial X(\gamma)^+\to \partial X(\gamma)$ as defined above is surjective.
\end{corollary}

\begin{proof}
By \Cref{Bowditch Prop 2.3.3} (2) and \hyperlink{thrmB}{Theorem B} (2), we have that $\partial X(\gamma) = \partial i_\gamma(\partial H)$. Suppose $y\in \partial X(\gamma)$. By \hyperlink{thrmB}{Theorem B} (2), there exists $a\in \partial H$ such that $\partial i_\gamma(a) = y$. Then by definition of $\tau^+$, we have that $\tau^+(\partial i_\gamma^+(a)) = \partial i_\gamma(a) = y$.
\end{proof}

\section{Ending Laminations}\label{EndingLaminations}

Recall that by \Cref{ShortExactSequenceConvention} we have fixed a short exact sequence $1\to H\to G\to Q\to 1$ of three infinite word-hyperbolic groups with Cayley graphs $\Gamma_H$, $\Gamma_G$, and $\Gamma_Q$, respectively. For each $g\in G$, conjugation by $g$ gives an automorphism $\phi_g$ of $H$ defined by $\phi_g(h) = g^{-1}hg$. Note that $\phi_g$ provides a bijection of the vertices of $\Gamma_H$ which is a quasi-isometry of $\Gamma_H$ with parameters depending on $|g|$. As such, $\phi_g$ extends to a homeomorphism of $\partial H$ that coincides with the action of left-multiplication by $g^{-1}$. We will also denote this homeomorphism by $\phi_g$. When $\lambda = [a,b]$ is a geodesic segment in $\Gamma_H$, we will denote a geodesic in $\Gamma_H$ between $\phi_g(a)$ and $\phi_g(b)$ by $\lambda_g$. Similarly, if $\lambda = (u,v)$ is a bi-infinite geodesic in $\Gamma_H$ with endpoints in $\partial H$, then $\lambda_g = (\phi_g(u), \phi_g(v)) = (g^{-1}u, g^{-1}v)$ also denotes the bi-infinite geodesic in $\Gamma_H$ between the images of the endpoints of $\lambda$ under the homeomorphism $\phi_g$.

Given $\kappa\geq 1$ and $\epsilon\geq 0$, define a \textit{$(\kappa, \epsilon)$-quasi-isometric section} to be a $(\kappa, \epsilon)$-quasi-isometric embedding $\sigma: \Gamma_Q\to \Gamma_G$ such that $P\cdot \sigma$ is the identity map on $\Gamma_Q$. The existence of such a quasi-isometric section in the setting of \Cref{ShortExactSequenceConvention} is guaranteed by Mosher \cite{Mos96}. If $\gamma \subseteq \Gamma_Q$ is a bi-infinite geodesic or a geodesic ray, we will also refer to a $(\kappa, \epsilon)$-quasi-isometric embedding $\sigma: \gamma\to \Gamma_G$ as a quasi-isometric section. All sections we consider in this paper are assumed to take vertices to vertices and edges to edge-paths. 

\begin{defn}
An \textit{algebraic lamination} on $H$ is defined to be a non-empty subset $L$ of the double boundary $\partial^2 H$ which is closed, symmetric (flip-invariant), and $H$-invariant. If $L\subseteq \partial^2 H$ is an algebraic lamination, an element $(p,q)\in L$ will be referred to as a \textit{leaf} of the lamination. As each point $(p,q)\in \partial^2 H$ can be represented by a bi-infinite geodesic $\lambda$ in $\Gamma_H$ from $p$ to $q$, we will sometimes refer to the geodesic $\lambda$ as a leaf of the lamination as well. 
\end{defn}

In \cite{Mit97}, Mitra describes a set of algebraic ending laminations on $\Gamma_H$ associated to the hyperbolic group extension (\ref{SES}) which are parametrized by points in the Gromov boundary of $\Gamma_Q$. These algebraic ending laminations are defined below.

\begin{convention}\label{Quasi-isometricSectionConvention}
Fix $\kappa\geq 1$ and $\epsilon\geq 0$, and let $\sigma: \Gamma_Q\to \Gamma_G$ be a quasi-isometric section of $\Gamma_Q$ into $\Gamma_G$. For a fixed $z\in \partial Q$, let $[1,z)\subseteq \Gamma_Q$ be a geodesic ray from the identity to $z$. Denote the $n^{\text{th}}$ vertex along $[1,z)$ by $z_n$, and set $g_n := \sigma(z_n)$. 
\end{convention}

\begin{defn}[Mitra, \cite{Mit97}]\label{LaminationDefinition} Let $z\in \partial Q$.
\begin{enumerate}
\item Let $h\in H$ be an element of infinite order. Choose a geodesic $[1,z)$ as in \Cref{Quasi-isometricSectionConvention}. 
Define $R_{z,h}$ to be the set of all pairs $(a,aw)\in H\times H$ such that there is some $n\geq 0$ for which $w\in [g_nhg_n^{-1}]_H$ and $w$ is a conjugacy minimal representative of $g_nhg_n^{-1}$ in $H$. Let $\overline{R_{z,h}}$ denote the closure of $R_{z,h}$ in $\widehat{H}\times \widehat{H}$, and set 
\[
\Lambda_{z,h}:= \overline{R_{z,h}} \cap \partial^2 H.
\]

So, $\Lambda_{z,h}$ consists of all points $(p,q)\in \partial^2 H$ for which there exists a sequence $(a_{n_i}, a_{n_i}w_{n_i})\in H\times H$ such that $(a_{n_i}, a_{n_i}w_{n_i})$ converges to $(p,q)$ in $\widehat{H}\times \widehat{H}$ as $n_i\to \infty$, where $w_{n_i}$ is some conjugacy minimal representative of $g_{n_i}hg_{n_i}^{-1}$ in $H$.

\item The \textit{algebraic ending lamination} corresponding to $z$ is 
\[
\Lambda_z:= \bigcup_{\substack{h\in H, \\ \text{$h$ infinite order}}} \Lambda_{z,h}.
\]
\item The \textit{algebraic ending lamination} for the short exact sequence (\ref{SES}) is
\[
\Lambda:= \bigcup_{z\in\partial Q} \Lambda_z.
\]
\end{enumerate}
\end{defn}

Note that $\Lambda_z$ is $H$-invariant and non-empty. While $\Lambda_{z,h}$ is not necessarily symmetric as defined, $\Lambda_{z,h}\cup \Lambda_{z,h^{-1}}$ is symmetric. Moreover, by \hyperlink{thrmC}{Theorem C} the subset $\Lambda_z\subseteq \partial^2 H$ is closed and therefore $\Lambda_z$ is an algebraic lamination on $H$. Mitra explained in \cite{Mit97} that in \Cref{LaminationDefinition} (2), it suffices to choose a finite collection of elements $h\in H$. Since $\Lambda_{z,h}$ is a closed subset of $\partial^2 H$ for each $h\in H$, this also shows that $\Lambda_z$ is closed.

\begin{remark}\label{LambdazIndependentOfRay}
We note the following about \Cref{LaminationDefinition} and the laminations $\Lambda_z$ and $\Lambda$.
\begin{enumerate}
    \item In \Cref{LaminationDefinition}, the quasi-isometric section $\sigma$ only needs to be defined on the ray $[1,z)$ rather than on all of $\Gamma_Q$.
    \item The lamination $\Lambda_z$ is independent of choice of quasi-isometric section, since if $\sigma: [1,z)$ and $\sigma': [1,z)$ are two quasi-isometric sections, $[\sigma(z_n)h\sigma(z_n)^{-1}]_H = [\sigma'(z_n)h\sigma'(z_n)^{-1}]_H$.
    \item The lamination $\Lambda_z$ is independent of geodesic ray $[1,z)$ by Mitra's Lemma 3.3 of \cite{Mit97}.
    \item The definitions of $\Lambda_z$ and $\Lambda$ are independent of the choice of generating set for $Q$. This follows from the proof of Lemma 3.3 \cite{Mit97} which can be adapted to show that $\Lambda_z$ is actually independent of quasigeodesic ray from 1 to $z$.
    \item Fix $z_0\in \Gamma_Q$, $z\in \partial Q$, and let $\gamma = [z_0,z)$ be a geodesic ray in $\Gamma_Q$ with vertices $z_n'\in \gamma$ such that $d_Q(z_0,z_n) = n$. Let $\sigma': [z_0,z) \to \Gamma_G$ be a quasi-isometric section with $\sigma'(z_n') = g_n'$ and let $\Lambda_z'$ be the algebraic ending lamination obtained by considering conjugacy minimal representatives of $g_n'h(g_n')^{-1}$. The proof of Lemma 3.3 \cite{Mit97} also shows that $\Lambda_z = \Lambda_z'$. So, when defining $\Lambda_z$, we can consider a geodesic ray from any basepoint $z_0\in \Gamma_Q$ converging to $z\in \partial Q$.

\end{enumerate}
\end{remark}

The next proposition shows how leaves of the lamination $\Lambda_z$ behave under the action of conjugation by elements of $G$.

\begin{proposition}\label{LambdaZ0}
Let $1\to H \to G \to Q \to 1$ be as in \Cref{ShortExactSequenceConvention} and let $P: \Gamma_G\to \Gamma_Q$ be the induced map. Then for all $g\in G$, $z\in \partial Q$, and $(u,v)\in\partial^2 H$, we have that $(u,v)$ is a leaf of $\Lambda_z$ if and only if $(g^{-1}u, g^{-1}v)$ is a leaf of $\Lambda_{P(g)^{-1}z}$.
\end{proposition}

\begin{proof}
Fix $z\in \partial Q$, $g\in G$, and set $q_0:= P(g)$. Let $\lambda = (u,v)$ be a leaf of $\Lambda_z$. If $[1,z)$ is a geodesic ray in $\Gamma_Q$ with vertices $1,z_1, z_2, \ldots$, then $q_0^{-1}\cdot [1,z) = [q_0^{-1}, q_0^{-1}z)$ is a geodesic ray in $\Gamma_Q$ with vertices $q_0^{-1}, q_0^{-1}z_1, q_0^{-1}z_2, \ldots$. Since $\Lambda_z$ is independent of quasi-isometric section, we may assume that $\sigma$ is a quasi-isometric section with $\sigma(q_0) = g$. As in \Cref{Quasi-isometricSectionConvention}, we will denote $\sigma(z_i)$ by $g_i$. 

Since $(u,v)\in \Lambda_z$, there is some sequence $(a_i, a_iw_i)\in H\times H$ such that $w_i\in [g_{n_i}hg_{n_i}^{-1}]_H$ is a conjugacy minimal representative of $g_{n_i}hg_{n_i}^{-1}$ in $H$ for some $n_i\geq 0$ and such that $a_i\to u$ and $a_iw_i\to v$ in $\widehat{\Gamma_H}$ as $i\to\infty$. Note that the sequence $(\phi_g(a_i), \phi_g(a_iw_i)) = (\phi_g(a_i), \phi_g(a_i)\phi_g(w_i))$ converges to $(\phi_g(u), \phi_g(v)) = (g^{-1}u, g^{-1}v)$ in $\widehat{H}\times \widehat{H}$. 

Since $w_i\in [g_{n_i}hg_{n_i}^{-1}]_H$, we have that $\phi_g(w_i)\in [ g^{-1}g_{n_i}hg_{n_i}^{-1}g ]_H$. As mentioned earlier, there exist constants $K\geq 1$ and $C\geq 0$ such that $\phi_{g}$ is a $(K, C)$-quasi-isometry. Since for each $i\geq 0$ we have that $w_i$ is a conjugacy minimal representative, \Cref{FHRQGeo} implies that there exists some $\kappa\geq 0$ such that for all $i\geq 0$, $\phi_{g}(w_i)$ is a $\kappa$-almost conjugacy minimal representative of $[ g^{-1}g_{n_i}hg_{n_i}^{-1}g ]_H$ in $H$. So, for each $i\geq 0$, there exists some $c_i\in H$ with $|c_i|_H\leq \kappa$ such that $c_i^{-1}\phi_g(w_i)c_i$ is a conjugacy minimal representative of $[g^{-1}g_{n_i}hg_{n_i}^{-1}g ]_H$. As $(\phi_g(a_i), \phi_g(a_i)\phi_g(w_i))\to (g^{-1}u, g^{-1}v)$ and $|c_i|\leq \kappa$ for all $i\geq 0$, we must also have that $(\phi_g(a_i)c_i, \phi_g(a_i)\phi_g(w_i)c_i) = (\phi_g(a_i)c_i, \phi_g(a_i)c_ic_i^{-1}\phi_g(w_i)c_i) \to (g^{-1}u,g^{-1}v)$. 

For each $n_i\geq 0$, the element $g^{-1}g_{n_i}$ is in the same coset of $H$ in $G$ as $\sigma(q_0^{-1}z_{n_i})$. So, $c_i^{-1}\phi_g(w_i)c_i$ is a conjugacy minimal representative of $[\sigma(q_0^{-1}z_{n_i}) h \sigma(q_0^{-1}z_{n_i})^{-1})]_H$. Therefore, by definition of $\Lambda_{q_0^{-1}z}$ and \Cref{LambdazIndependentOfRay} (5), we have that $\lambda_g = (g^{-1}u, g^{-1}v)$ is a leaf of $\Lambda_{q_0^{-1}z}$.

Now, suppose that $\lambda_g = (g^{-1}u,g^{-1}n)$ is a leaf of $\Lambda_{P(g)^{-1}z}$. Let $g^{-1}u = u'$, $g^{-1}v = v'$, and let $\lambda' = (u',v')$. Then the forward direction of this proposition shows that $\lambda'_{g^{-1}}\in \Lambda_{P(g^{-1})^{-1}P(g)^{-1}z} = \Lambda_z$. As $\lambda'_{g^{-1}} = (u,v) = \lambda$, the reverse direction of this proposition follows.
\end{proof}

The main result of Mitra in \cite{Mit97} is the following.
\begin{theorem}[Mitra \cite{Mit97}, Theorem 4.11]\label{Mitra, Theorem 4.11}
Suppose that $1\to H\to G\to Q\to 1$ is as in \Cref{ShortExactSequenceConvention} and $\partial i: \partial H \to \partial G$ is the Cannon-Thurston map. Then for distinct points $u,v\in \partial H$, $\partial i(u) = \partial i(v)$ if and only if $(u,v)\in \Lambda$. 
\end{theorem}

The goal of the remainder of this section is to prove \hyperlink{thrmC}{Theorem C}. We first show that if $\lambda = (u,v)$ is a leaf of $\Lambda_z$, then $\partial i_\gamma^+$ identifies the endpoints $u$ and $v$.

\begin{proposition}\label{Injectivity}
Let $1\to H\to G\to Q\to 1$ be as in \Cref{ShortExactSequenceConvention}, $\gamma$ be as in \Cref{StackConvention}, $i_\gamma^+$ be as in \Cref{StackInclusionConvention}, and let $\partial i_\gamma^+:\partial H \rightarrow \partial X(\gamma)^+$ denote the Cannon-Thurston map. If $\lambda = (u,v)$ is a leaf of $\Lambda_z$, then $\partial i_\gamma^+(u) = \partial i_\gamma^+(v)$.
\end{proposition}

\begin{proof}
Let $\lambda = (u,v)\in \Lambda_z$ and suppose that $h\in H$ is such that $\lambda$ is a leaf of $\Lambda_{z, h}$. By \Cref{LambdazIndependentOfRay}, we can consider $\Lambda_{z,h}$ defined by the geodesic ray $[z_0,z)$. If $\sigma': \Gamma_Q\to \Gamma_G$ is any quasi-isometric section, then $[\sigma'(z_i)h\sigma'(z_i)^{-1}]_H = [g_ihg_i^{-1}]_H$ for all $z_i\in [z_0,z)$. Hence, there exist elements $a_i\in H$ and conjugacy minimal representatives $w_i\in [g_{n_i}hg_{n_i}^{-1}]_H$ for some $n_i\geq 0$ such that $a_i\to u$ and $a_iw_i\to v$ as $i\to \infty$. Note that since $w_i$ is conjugacy minimal, we have that $[a_iw_i^{-1}, a_i]\cup [a_i, a_iw_i]\cup [a_iw_i, a_iw_i^2]$ is a $(1,C_1)$-quasigeodesic for $C_1 = C_1(\delta)$ by \Cref{ConjMinBoundedGromov} and \Cref{SmallGromovProductsQuasigeodesic}. So, we have that $a_iw_i^{-1}\to u$ and $a_iw_i^2\to v$ as well.

Suppose for each $i\geq 0$ that $g_{n_i}hg_{n_i}^{-1} =_H c_i^{-1}w_i'c_i$, where $c_i\in H$ is a minimal length element conjugating $g_{n_i}hg_{n_i}^{-1}$ to a cyclic conjugate $w_i'$ of $w_i$. Mark vertices $p_i$ on $[a_iw_i^{-1}, a_i]$ and $q_i$ on $[a_iw_i, a_iw_i^2]$ where the path labeled by $(w_i')^2$ begins and ends. Let $x_i = p_ic_i$ and $y_i = q_ic_i$ denote the vertices at the end of the paths labeled by $c_i$ which start at $p_i$ and $q_i$, respectively. Now, as in the proof of \Cref{GeneralWordAsFHR} the minimality of $|c_i|$ requires that $(x_i, q_i; \Gamma_H)_{p_i} \leq \delta$ and $(p_i,y_i; \Gamma_H)_{q_i}\leq \delta$. Hence, by \Cref{SmallGromovProductsQuasigeodesic}, $[x_i, p_i]\cup [p_i, q_i] \cup [q_i, y_i]$ is a $(1,8\delta)$-quasigeodesic in $\Gamma_H$. So, we must have that $x_i\to u$ and $y_i\to v$ in $\widehat{\Gamma_H}$. Note that the geodesic in $\Gamma_H$ between $x_i$ and $y_i$ is labeled by the word $c_i^{-1}(w_i')^2 c_i =_H g_{n_i}h^2g_{n_i}^{-1}$.

Recall that $i_\gamma^+(x_i) = x_ig_0$ and $i_\gamma^+(y_i) = y_ig_0$. So, the geodesic between $i_\gamma^+(x_i)$ and $i_\gamma^+(y_i)$ in $X(\gamma)$ is labeled by a word representing the element $g_0^{-1}g_{n_i}h^2g_{n_i}^{-1}g_0$. To show that $\partial i_\gamma^+(u) = \partial i_\gamma^+(v)$, we must show that in $X(\gamma)^+$, the distance between some fixed point and the geodesic between $i_\gamma^+(x_{i})$ and $i_\gamma^+(y_{i})$, goes to infinity as $i \to \infty$. For each $i > 0$, consider the path $\rho_{i} \in X(\gamma)^+$ from $i_\gamma^+(x_{i})$ to $i_\gamma^+(y_{i})$ 
which consists of the geodesic $[x_ig_0, x_ig_{n_i}]$ labeled by $g_0^{-1}g_{n_i}$, followed by the quasigeodesic from $x_ig_{n_i}$ to $x_ig_{n_i}h^2$ labeled by $h^2$, followed by the geodesic $[x_ig_{n_i}h^2, y_ig_0]$ labeled by $g_{n_i}^{-1}g_0$. This path is a $(1,C)$-quasigeodesic in $X(\gamma)^+$ by \Cref{QuasigeodesicLemma} for some constant $C\geq 0$ independent of $i$.

So, take an arbitrary point $p\in \rho_{n}$. We will show that $p$ is far from $g_0$ in $X(\gamma)^+$, and so the distance in $X(\gamma)^+$ between a quasigeodesic between $i_\gamma^+(x_n)$ and $i_\gamma^+(y_n)$ and $g_0$ goes to infinity as $n$ goes to infinity. Note that since $d_H(1,x_n)\to \infty$, we must have that $d_{X(\gamma)^+}(g_0, i_\gamma^+(x_n))\to \infty$. 

Suppose first that the point $p$ belongs to the initial part of $\rho_{n}$ which is labeled by $g_0^{-1}g_{n}$. In this case, $p = i_\gamma^+(x_n)g_0^{-1}g_j$, $0\leq j\leq n$. There are two cases for us to consider:
\begin{enumerate}
\item If $j \leq \frac{1}{2} d_{X(\gamma)^+}(g_0, i_\gamma^+(x_n))$, then\begin{align*}
d_{X(\gamma)^+}(p, g_0) &= d_{X(\gamma)^+}(i_\gamma^+(x_n)g_0^{-1}g_j, g_0)\\
&\geq d_{X(\gamma)^+}(i_\gamma^+(x_n), g_0) - j \\
&\geq  \frac{1}{2}d_{X(\gamma)^+}(g_0, i_\gamma^+(x_n)). \end{align*}
\item If $j > \frac{1}{2}d_{X(\gamma)^+}(g_0, i_\gamma^+(x_n))$, then \begin{align*} d_{X(\gamma)^+}(p,g_0) &= d_{X(\gamma)^+}(i_\gamma^+(x_n) g_j, g_0)\\
&\geq j  \\ &> \frac{1}{2}d_{X(\gamma)^+}(g_0, i_\gamma^+(x_n)).
\end{align*}
\end{enumerate} In both cases, $d_{X(\gamma)^+}(p, g_0)\to \infty$ as $n\to \infty$. The case where $p$ belongs to the terminal part of $\rho_n$ which is labeled by $g_{n}^{-1}g_0$ is handled similarly.

Finally, if $p$ is a vertex in the portion of $\rho_n$ which is labeled by $h^2$, then since $h^2\in H$ is fixed, in $X(\gamma)^+$, $p$ must lie a bounded distance away from the element $i_\gamma^+(x_n)g_0^{-1}g_n$. In this case, we have that $d_{X(\gamma)^+}(p,g_0) \geq d_{X(\gamma)^+}(i_\gamma^+(x_n)g_0^{-1}g_n, g_0) - |h^2|_{H}$, and $d_{X(\gamma)^+}(i_\gamma^+(x_n)g_0^{-1}g_n, g_0) - |h^2|_{H}\to \infty$ as $n\to \infty$.  Therefore, the distance between $[i_\gamma(x_n)^+, i_\gamma^+(y_n)]_{X(\gamma)^+}$ and $g_0$ in $X(\gamma)^+$ goes to infinity as $n\to \infty$. Hence, $\partial i_\gamma^+(u) = \partial i_\gamma^+(v)$.
\end{proof}

The following several lemmas from Mitra \cite{Mit97} will allow us to show that certain geodesics are conjugacy minimal representatives. We have stated and proved these results in the setting where $\gamma = (z',z)\subseteq \Gamma_Q$ does not necessarily go through the identity. However, we will apply these results in a simpler setting where $\gamma$ does go through the identity. We have included the more general statements here to illuminate what happens in the general setting. The following three results are used to prove \Cref{Analogue of 4.5}, which is key to the proof of \hyperlink{thrmC}{Theorem C}.

\begin{lemma}[Cf Mitra \cite{Mit97}, Lemma 4.2 \cite{Mit97}]\label{Analogue of 4.2}
There exists $\kappa \geq 0$ such that for any $(u,v)\in \partial^2 H$ with $\partial i(u) = \partial i(v)$, any geodesic subsegment $[p,q]$ of $\lambda = (u,v)$ has an extension $[r,q]$ in $\lambda$ with $d_H(p,r)$ equal to 0 or 1 such that $[r,q]$ is a $\kappa$-almost conjugacy minimal representative.
\end{lemma}

The next lemma is proved in a similar manner to Mitra's Lemma 4.3 in \cite{Mit97}.

\begin{lemma}\label{Analogue of 4.3}
Given $\kappa\geq 0$, there exists $C\geq 1$ such that for any distinct $z,z'\in \partial Q$ and for any geodesic $\gamma = (z', z)\subset \Gamma_Q$ with $z_0\in \gamma$ the following holds:

If $\lambda = [1,h]\subseteq \Gamma_H$ and $\lambda_{g_0}$ is a $\kappa$-almost conjugacy minimal representative for some $g_0\in P^{-1}(z_0)$, then there exists a $(C, 0)$-quasi-isometric section $\sigma_0$ of $(z', z)$ into $X(\gamma)$ containing $g_0$ such that for all $g\neq g_0$ in $\sigma_0((z', z))$, $\lambda_{g}$ is a conjugacy minimal representative.
\end{lemma}

\begin{proof}
Let $\gamma = (z',z)$ be as in \Cref{StackConvention}. Let $\sigma: (z',z) \to X(\gamma)$ be an isometric lift of $(z',z)$ into $X(\gamma)$ with $\sigma(z_0) = g_0$ and such that $\lambda_{g_0}$ is a $\kappa$-almost conjugacy minimal representative for some $\kappa \geq 0$. We will construct the quasi-isometric section $\sigma_0$ satisfying the conclusions of the lemma inductively. 

Set $\sigma_0(z_0) = g_0$. For each $n\geq 0$ set $s_n:= \sigma(z_n)^{-1}\sigma(z_{n+1})$, and for each $n\leq 0$ set $s_{n-1} = \sigma(z_n)^{-1}\sigma(z_{n-1})$. Note that since $\sigma$ is an isometric embedding, $|s_n| = 1$ for all $n$. So, there exists some $K_1\geq 1$ and $\epsilon_1\geq 0$ such that $\phi_{s_n}: \Gamma_H\to \Gamma_H$ is a $(K_1, \epsilon_1)$-quasi-isometry for all $n\geq 0$. As $\lambda_{g_0}$ is a $\kappa$-almost conjugacy minimal representative, there exists $\kappa'\geq 0$ such that $\phi_{s_0}(\lambda_{g_0}) = \lambda_{g_0s_0}$ is a $\kappa'$-conjugacy minimal representative by \Cref{FHRQGeo}. By \Cref{ShortElementToConjMin}, there exists $c_0\in H$ and $M'\geq 0$ with $|c_0|_H\leq M'$ such that $\lambda_{g_0s_0c_0}$ is a conjugacy minimal representative. Set $\sigma_0(z_1):= g_0s_0c_0$. We can similarly define $\sigma_0(z_{-1})$.

Suppose that $\sigma_0(z_j)$ has been constructed satisfying the conclusions of the lemma for all $-m\leq j\leq n$. By assumption, $\lambda_{\sigma_0(z_n)}$ is a conjugacy minimal representative, and so by \Cref{FHRQGeo} there exists $\kappa''\geq 0$ such that $\lambda_{\sigma_0(z_n)s_n}$ is a $\kappa''$-almost conjugacy minimal representative. Then by \Cref{ShortElementToConjMin}, there exists $c_n\in H$ and $M''\geq 0$ with $|c_n|_H\leq M''$ such that $\lambda_{\sigma_0(z_n)s_nc_n}$ is a conjugacy minimal representative. Set $\sigma_0(z_{n+1}):= \sigma_0(z_n)s_nc_n$. We can similarly define $\sigma_0(z_{-m-1})$. Note that $d_{X(\gamma)}(\sigma_0(z_i), \sigma_0(z_{i+1})) \leq \max\{M',M''\}$, and so $\sigma_0$ is a $(C,0)$-quasi-isometric section, where $C := \max\{M',M''\}$ and $\lambda_g$ is a conjugacy minimal representative for all $g\neq g_0$ in $\sigma_0((z',z))$.
\end{proof}

The following corollary is obtained from the previous lemma by translating the quasi-isometric section by an element of $G$. Here, we choose the quasi-isometric section $\sigma_0$ to go through the point $g_0\in \Gamma_G$ rather than the identity.

\begin{corollary}[Cf Mitra \cite{Mit97} Corollary 4.4]\label{Analogue of 4.4}
Given $\kappa \geq 0$, there exists $C\geq 1$ such that for any geodesic ray $[z_0, z)$ in $\Gamma_Q$ and any $g\in P^{-1}([z_0,z))$ the following holds:

If $\lambda = [1,h]\subseteq \Gamma_H$ and $\lambda_{g_0}$ is a $\kappa$-almost conjugacy minimal representative for some $g_0\in P^{-1}(z_0)$, then there exists a $(C,0)$-quasi-isometric section $\sigma_0$ of $[z_0, z)$ into $\Gamma_G$ containing $g\in \Gamma_G$ such that for all $g'\neq g$ in $\sigma_0([z_0, z))$, $\lambda_{g_0g^{-1}g'}$ is a conjugacy minimal representative.
\end{corollary}

\begin{proof}
By \Cref{Analogue of 4.3}, there exists a $(C,0)$-quasi-isometric section $\sigma': (z',z)\to X(\gamma)$ with $\sigma'(z_0) = g_0$ such that for all $g'\neq g_0$ in $\sigma'((z',z))$, $\lambda_{g'}$ is a conjugacy minimal representative. Suppose that $g\in P^{-1}(z_n)$ and set $\sigma_0(z_n) := g$. For each integer $i$ with $i\geq -n$, set $\sigma_0(z_{n+i}) := t_{gg_0^{-1}}\cdot \sigma'(z_i)$. Now, $\sigma_0: [z_0,z)\to X(\gamma)^+$ is a $(C,0)$-quasi-isometric section since it is a left-translate of $\sigma'$ by $gg_0^{-1}\in G$. Also, note that for all $g'\neq g$ in $\sigma([z_0,z))$, we have that $g' = t_{gg_0^{-1}}\cdot \sigma'(z_i)$ for some $i\geq -n$ with $i\neq 0$. Then, $\lambda_{g_0g^{-1}g'} = \lambda_{g_0g^{-1}gg_0^{-1}\sigma'(z_{i})} = \lambda_{\sigma'(z_{i})}$ is a conjugacy minimal representative by \Cref{Analogue of 4.3}.
\end{proof}

The following lemma will allow us to reduce to the simpler setting where $\gamma = (z',z)\subseteq \Gamma_Q$ passes through the identity in $\Gamma_Q$.

\begin{lemma}\label{FormalReduction}
Suppose $\gamma = (z',z)$ is as in \Cref{StackConvention} and let $\gamma' := z_0^{-1} \cdot \gamma = (z_0^{-1}z', z_0^{-1}z)$ be as in \Cref{IdentityStackConvention}. Let $X(\gamma)$ and $X(\gamma')$ be the stacks as in \Cref{StackInclusionConvention} where the section $\sigma: \gamma\to X(\gamma)$ is such that $\sigma(z_0) = g_0$ and $\sigma': \gamma'\to X(\gamma')$ is chosen so that $\sigma'= g_0^{-1}\cdot \sigma$. Let $i_\gamma^+$ and $i_{\gamma'}^+$ be as in \Cref{StackInclusionConvention}, and let $\partial i_\gamma^+: \partial H\to \partial X(\gamma)^+$ and $\partial i_{\gamma'}^+: \partial H\to \partial X(\gamma')^+$ be the Cannon-Thurston maps. 

Then for any two distinct points $u,v\in \partial H$, $\partial i_\gamma^+(u) = \partial i_\gamma^+(v)$ if and only if $\partial i_{\gamma'}^+(\phi_{g_0}(u)) = \partial i_{\gamma'}^+(\phi_{g_0}(v))$, where $g_0\in P^{-1}(z_0)$.
\end{lemma}

\begin{proof}
Let $\gamma = (z',z)$ be as in \Cref{StackConvention}, let $\gamma' := z_0^{-1} \cdot \gamma = (z_0^{-1}z', z_0^{-1}z)$, and fix some $g_0\in P^{-1}(z_0)$. Recall that $i_\gamma^+$ is given by $i_\gamma^+(h) = t_{g_0} \cdot \phi_{g_0}(h) = hg_0$ and $i_{\gamma'}^+$ is given by $i_{\gamma'}(h) = h$. Suppose first that $u,v\in \partial H$ are distinct points such that $\partial i_\gamma^+(u) = \partial i_\gamma^+(v)$. Then, for any sequences $(u_n), (v_n)\in \Gamma_H$ with $u_n\to u$ and $v_n\to v$ in $\widehat{\Gamma_H}$, we have that in $\widehat{X(\gamma)^+}$, $\lim_{n\to \infty} i_\gamma^+(u_n) = \lim_{i\to \infty} i_\gamma^+(v_n)$. So, in $\widehat{X(\gamma)^+}$ we have that $\lim_{i\to \infty} u_ng_0 = \lim_{i\to \infty} v_ng_0$. Note that $X(\gamma)^+ = g_0X(\gamma')^+$ and so left-translation by $g_0^{-1}$ gives an isometry from $X(\gamma)^+$ to $X(\gamma')^+$. Therefore, in $\widehat{X(\gamma')^+}$ we have that $\lim_{i\to \infty} g_0^{-1}u_ng_0 = \lim_{i\to \infty} g_0^{-1}v_ng_0$. So by definition of $i_{\gamma'}^+$, we have that $\lim_{i\to \infty} i_{\gamma'}^+(\phi_{g_0} (u_n)) = \lim_{i\to \infty} i_{\gamma'}^+(\phi_{g_0} (v_n))$ in $\widehat{X(\gamma')^+}$. Since in $\widehat{\Gamma_H}$ $\phi_{g_0}(u_n)\to \phi_{g_0}(u)$ and $\phi_{g_0}(v_n) \to \phi_{g_0}(v)$ as $n\to \infty$, we have that $\partial i_{\gamma'}^+(\phi_{g_0}(u)) = \partial i_{\gamma'}^+(\phi_{g_0}(v))$ by the continuity of $\widehat{i}_{\gamma'}^+$ (\Cref{i_gamma Cannon-Thurston}). The reverse implication follows in the same manner by noting that left-translation by $g_0$ gives an isometry from $X(\gamma')^+$ to $X(\gamma)^+$.
\end{proof}

The following result follows directly from \Cref{Analogue of 4.2} and \Cref{Analogue of 4.4}. This corollary will be used in the proof of \hyperlink{thrmC}{Theorem C} to construct a sequence of conjugacy minimal representatives which converge to some bi-infinite geodesic $\lambda\subseteq \partial^2 H$ whose endpoints are identified by $\partial i_\gamma^+$.

\begin{corollary}[Cf Mitra \cite{Mit97} Lemma 4.5]\label{Analogue of 4.5}
There exists $C'$ such that for any $\lambda = (u,v)$, $u,v\in \partial H$ with $\partial i_\gamma^+(u) = \partial i_\gamma^+(v)$, any geodesic ray $[z_0, z)$ in $\Gamma_Q$, and any geodesic subsegment $[p,q]$ of $\lambda_g$ for some $g\in P^{-1}([z_0, z))$ the following holds:

There exists an extension $[r,q] = \mu$ of $[p,q]$ in $\lambda_g$ with $d_H(p,r)$ equal to 0 or 1 and a $(C', 0)$-quasi-isometric section $\sigma: [z_0, z)\to X(\gamma)$ such that $gr\in \sigma([z_0,z))$ and $\mu_{g_0r^{-1}g^{-1}g'}$ is a conjugacy minimal representative for all $g'\neq gr$ in $\sigma([z_0, z))$.
\end{corollary}

\begin{proof}
Let $\lambda = (u,v)$ be such that $\partial i_\gamma^+(u) = \partial i_\gamma^+(v)$, let $[z_0,z)\in \Gamma_Q$ be a geodesic ray, let $g\in P^{-1}([z_0,z))$, and let $[p,q]$ be any geodesic subsegment of $\lambda_g = (\phi_g(u), \phi_g(v))$. 
By \Cref{FormalReduction}, $\partial i_{\gamma'}^+(\phi_{g_0}(u)) = \partial i_{\gamma'}^+(\phi_{g_0}(v))$. So by \Cref{igamma=i}, we have that $\partial i(\phi_g(u)) = \partial i(\phi_g(v))$. 
So by \Cref{Analogue of 4.2}, there exists an extension $[r,q] = \mu$ of $[p,q]$ in $\lambda_g$ with $d_H(p,r)$ equal to 0 or 1 and such that $[r,q]$ is a $\kappa$-almost conjugacy minimal representative for some $\kappa \geq 0$. Let $\mu' = [1, r^{-1}q]$ and note that $\mu'$ is also a $\kappa$-almost conjugacy minimal representative since it has the same label as $\mu$. By \Cref{FHRQGeo}, $\mu_{g_0}$ and $\mu'_{g_0}$ are $\kappa'$-almost conjugacy minimal representatives for some $\kappa'\geq 0$ depending on $g_0$. So by \Cref{Analogue of 4.4}, there exists $C'\geq 1$ and a $(C', 0)$-quasi-isometric section $\sigma: [z_0,z)\to \Gamma_G$ containing $gr\in \Gamma_G$ such that for all $g'\neq gr\in \sigma([z_0,z))$, $\mu'_{g_0r^{-1}g^{-1}g'}$ is a conjugacy minimal representative. Therefore, $\mu_{g_0r^{-1}g^{-1}g'}$ is also a conjugacy minimal representative.  
\end{proof}

For the next portion of this section, we will assume that the bi-infinite geodesic $\gamma = (z',z)\subseteq \Gamma_Q$ goes through the identity in $Q$, and so $\gamma^+ = [1,z)$. Note that several of the previous lemmas simplify in this case. We now make the following convention.

\begin{convention}\label{FormalConventionThroughIdentity}
Let $\gamma = (z',z)$ be a bi-infinite geodesic in $\Gamma_Q$ between $z',z\in \partial Q$ with $z'\neq z$ and assume that $1\in \gamma$. Label the sequence of vertices in order along the portion of $\gamma$ from 1 to $z$ by $1 = z_0, z_1, z_2, \ldots$. Similarly, label the sequence of vertices in order along the portion of $\gamma$ from 1 to $z'$ by $1 = z_0, z_{-1}, z_{-2}, \ldots$. Let $\sigma_0: \gamma\to \Gamma_G$ denote an isometric lift of $\gamma$ through the identity in $\Gamma_G$, i.e. such that $\sigma_0(1) = 1$, and set $g_i:= \sigma_0(z_i)$. Let $X(\gamma)$ and $X(\gamma)^+$ denote the stacks over $\gamma = (z',z)$ and $\gamma^+ = [1,z)$, respectively. Finally, let $i_\gamma: \Gamma_H\to X(\gamma)$ and $i_\gamma^+: \Gamma_H\to X(\gamma)^+$ be the respective inclusion maps given by $i_\gamma(h) = h$ and $i_\gamma^+(h) = h$ for all $h\in H$. 
\end{convention}

Before proving \hyperlink{thrmC}{Theorem C}, we will first introduce some necessary terminology as well as some lemmas which were first stated by Mitra in \cite{Mit97}.

Given a (finite or infinite) geodesic $\lambda\subset \widehat{\Gamma_H}$ with endpoints $a,b\in \widehat{\Gamma_H}$ and an element $g\in G$, recall that $\lambda_g\subset \widehat{\Gamma_H}$ denotes the geodesic joining $\phi_g(a) = g^{-1}ag$ and $\phi_g(b) = g^{-1}bg$. For any quasi-isometric section $\sigma: \Gamma_Q \to \Gamma_G$ and geodesic $\lambda$, Mitra defines the set 
$$B(\lambda, \sigma) := \bigcup_{g\in \sigma(Q)} t_g\cdot i(\lambda_g),$$ 
where $t_g$ denotes left-translation by the element $g\in G$. For our purposes, we will consider the subset of $B(\lambda, \sigma)$ which lives in $X(\gamma)^+$: $$B_{\gamma^+}(\lambda, \sigma) = \bigcup_{g\in \sigma([1, z))} t_g\cdot i_\gamma^+(\lambda_g).$$ 
Note that $B_{\gamma^+}(\lambda,\sigma) = B(\lambda, \sigma)\cap P^{-1}([1, z))$ and that if $\lambda$ is a bi-infinite geodesic, then $B_{\gamma^+}(\lambda,\sigma)$ is independent of quasi-isometric section $\sigma$ for the same reason Mitra uses to show $B(\lambda,\sigma)$ is independent of quasi-isometric section \cite{Mit97}.

On the vertices of $\Gamma_H$, define the map $\pi_{g,\lambda}: \Gamma_H\to \lambda_g$ by sending $h\in H$ to a closest vertex on $\lambda_g$. We will now define a projection map to the set $B_{\gamma^+}(\lambda, \sigma)$. As $\sigma$ is a quasi-isometric section, for each $g'\in X(\gamma)^+$, there is a unique $g\in \sigma([1,z))$ and $h\in H$ such that $g'= t_{g}\cdot i_\gamma^+(h)$. So, define $$\Pi_\lambda^\sigma(g') = \Pi_\lambda^\sigma\cdot t_g\cdot i_\gamma^+(h) := t_g\cdot i_\gamma^+\cdot \pi_{g,\lambda}(h).$$

The following statements are versions of the analogous statements from Mitra \cite{Mit97} which apply to the setting in which we are working. In most cases, the proofs that Mitra provided go through with no changes to the reasoning. We provide details of the necessary modifications where they are needed.

The same proof of Mitra's Theorem 3.7 of \cite{Mit98} verifies the following statement. In particular, this lemma will be used to show that if $\sigma: [1,z)\to X(\gamma)^+$ is a $(K, \epsilon)$-quasi-isometric section, then the projection of $\sigma$ to $B_{\gamma^+}(\lambda, \sigma)$ is also a quasi-isometric section.

\begin{lemma}[Cf Mitra \cite{Mit97}, Theorem 4.6]\label{Analogue of 4.6}
For all $K\geq 1$ and $\epsilon\geq 0$, there exists a constant $C\geq 1$ such that if $\sigma: [1,z)\to X(\gamma)^+$ is any $(K,\epsilon)$-quasi-isometric section and $\lambda\subseteq \Gamma_H$ is any bi-infinite geodesic, then for all $x,y\in X(\gamma)^+$, $d_{X(\gamma)^+}(\Pi_\lambda^\sigma(x), \Pi_\lambda^\sigma(y)) \leq Cd_{X(\gamma)^+}(x,y)$.
\end{lemma}

\begin{lemma}[Cf Mitra \cite{Mit97}, Lemma 4.7]\label{Analogue of 4.7}
For all $K\geq 1$ and $\epsilon\geq 0$ there exists $A\geq 1$ such that if $\sigma: [1,z)\to X(\gamma)^+$ is a $(K,\epsilon)$-quasi-isometric section, then for all $p,q\in \sigma([1,z))$ and $x\in t_p\cdot i_\gamma^+(\lambda_p)$ there exists $y\in t_q\cdot i_\gamma^+(\lambda_q)$ such that $d_{X(\gamma)^+}(x,y)\leq Ad_Q(Px, Py) = Ad_Q(Pp, Pq)$. 
\end{lemma}

\begin{proof}
Let $\sigma: [1,z)\to X(\gamma)^+$ be a $(K,\epsilon)$-quasi-isometric section, $p,q\in \sigma([1,z))$, $x\in t_p\cdot i_\gamma^+(\lambda_p)$, and set $y = \Pi_\lambda^\sigma(xp^{-1}q)$. Note that $y\in t_q\cdot i_\gamma^+(\lambda_q)$. Then by \Cref{Analogue of 4.6}, there exists a constant $C\geq 1$ such that $d_{X(\gamma)^+}(\Pi_\lambda^\sigma(x), \Pi_\lambda^\sigma(xp^{-1}q)) = d_{X(\gamma)^+}(x,y) \leq Cd_{X(\gamma)^+}(x,xp^{-1}q)$. Since $p,q\in \sigma([1,z))$ and $\sigma$ is a $(K, \epsilon)$-quasi-isometric section, we have that $|p^{-1}q| \leq Kd_Q(Pp, Pq) + \epsilon$. Therefore, $d_{X(\gamma)^+}(x,xp^{-1}q) = |p^{-1}q| \leq Kd_Q(Pp, Pq) + \epsilon$. So, let $A = C(K + \epsilon)$. As $Px = Pp$ and $Py = Pq$, we have finally that $d_{X(\gamma)^+}(x,y)\leq A d_Q(Px,Py) = Ad_Q(Pp, Pq)$ as required.
\end{proof}

The following is the version of Lemma 4.8 \cite{Mit97} that we need for our purposes. It is proved by an argument similar to the one given by Mitra using the previous lemma.
\begin{lemma}[Cf Mitra \cite{Mit97}, Lemma 4.8]\label{Analogue of 4.8}
For all $K\geq 1$ and $\epsilon\geq 0$ there exists $M\geq 0$ such that the following holds. Suppose $\lambda$ is a bi-infinite geodesic in $\Gamma_H$ and $a$ is a vertex on $\lambda$ splitting $\lambda$ into semi-infinite geodesics $\lambda^-$ and $\lambda^+$. Suppose further that $\sigma: [1,z)\to X(\gamma)^+$ is a $(K, \epsilon)$-quasi-isometric section such that $\sigma([1,z))\subseteq B_{\gamma^+}(\lambda, \sigma)$ and $i_\gamma^+(a)\in \sigma([1,z))$.
Then, any geodesic in $X(\gamma)^+$ joining a point in $B_{\gamma^+}(\lambda^-, \sigma)$ to a point in $B_{\gamma^+}(\lambda^+, \sigma)$ passes through an $M$-neighborhood of $\sigma([1,z))$. 
\end{lemma}

\begin{lemma}[Cf Mitra \cite{Mit97}, Corollary 4.10]\label{Analogue of 4.10}
Given $K\geq 1$, $\epsilon\geq 0$, there exists $\alpha$ such that if $\lambda = (u,v)$ is such that $\partial i_\gamma^+(u) = \partial i_\gamma^+(v)$ then the following is satisfied: 

If $\sigma$ and $\sigma'$ are $(K, \epsilon)$-quasi-isometric sections such that $B_{\gamma^+}(\lambda, \sigma) = B_{\gamma^+}(\lambda, \sigma')$ and $\sigma, \sigma'$ are contained in $B_{\gamma^+}(\lambda, \sigma)$, then there exists $N\geq 0$ such that for all $n\geq N$, $$d_{X(\gamma)^+}(\sigma(z_n), \sigma'(z_n))\leq \alpha.$$
\end{lemma}
 
\begin{proof}
Let $\lambda = (u,v)$ be such that $\partial i_\gamma^+(u) = \partial i_\gamma^+(v)$ and let $\sigma$ and $\sigma'$ be $(K,\epsilon)$-quasi-isometric sections satisfying the hypotheses of the lemma. Let $(p_n)$ and $(q_n)$ be a sequence of vertices on $\lambda$ such that $p_n\to u$ and $q_n\to v$ as $n\to \infty$. For each $n\geq 0$, \Cref{Analogue of 4.8} guarantees there exist points $z_{n'}, z_{n''}\in [1,z)$ such that any geodesic in $X(\gamma)^+$ joining $i_\gamma^+(p_n)$ to $i_\gamma^+(q_n)$ passes through an $M$-neighborhood of both $\sigma(z_{n'})$ and $\sigma'(z_{n''})$. Since $\partial i_\gamma^+(u) = \partial i_\gamma^+(v)$ and $\widehat{i_\gamma^+}$ is continuous, we must have that the sequences $\{i_\gamma^+(p_n)\}$, $\{i_\gamma^+(q_n)\}$, $\{\sigma(z_{n'})\}$, and $\{\sigma'(z_{n''})\}$ all converge to the same point in $\partial X(\gamma)^+$. Since $\sigma$ and $\sigma'$ are quasi-isometric sections of $[1, z)$ into $X(\gamma)^+$ and as
$d_{X(\gamma)^+}(1,[i_\gamma^+(p_n), i_\gamma^+(q_n)])\to \infty$, we must have that $z_{n'}\to z$ and $z_{n''}\to z$. Therefore, $\sigma([1,z))$ and $\sigma'([1,z))$ are asymptotic quasigeodesic rays in $X(\gamma)^+$ and we have that for all $n\geq N$,
$$\max\{d_{X(\gamma)^+}(\sigma(z_n), \sigma'([1,z))), d_{X(\gamma)^+}(\sigma([1,z)), \sigma'(z_n)) \}\leq \alpha'.$$ 
But since $\sigma$ and $\sigma'$ are $(K,\epsilon)$-quasi-isometric sections, if $z_{n'}$ is such that $d_{X(\gamma)^+}(\sigma(z_n), \sigma'(z_{n'}))\leq \alpha'$, then we have that \begin{align*}
    d_{X(\gamma)^+}(\sigma(z_n), \sigma'(z_{n}))&\leq d_{X(\gamma)^+}(\sigma(z_n), \sigma'(z_{n'})) + d_{X(\gamma)^+}(\sigma'(z_{n'}), \sigma'(z_{n})) \\
    &\leq \alpha' + K|n - n'| + \epsilon \\
    &\leq \alpha' + K\alpha' + \epsilon = \alpha
\end{align*} 
Thus for all $n\geq N$, $d_{X(\gamma)^+}(\sigma(z_n), \sigma'(z_n))\leq \alpha.$
\end{proof}

We are now ready to prove the main theorem of this section which is reminiscent of Mitra's Theorem 4.11 \cite{Mit97}. 

\begin{mainthrmC}\hypertarget{thrmC}{}
Let $1\to H\to G\to Q\to 1$ be a short exact sequence of infinite, finitely generated, word-hyperbolic groups. Let $z,z'\in \partial Q$ be distinct and let $\gamma\subseteq \Gamma_Q$ be a bi-infinite geodesic in $\Gamma_Q$ between $z$ and $z'$. Let $i_\gamma^+: \Gamma_H\to X(\gamma)^+$ be the inclusion of $\Gamma_H$ into the semi-infinite stack $X(\gamma)^+$ over $\gamma^+ = [z_0,z)$ as in \Cref{StackInclusionConvention}, and let $\partial i_\gamma^+: \partial H\to \partial X(\gamma)^+$ be the Cannon-Thurston map. 

Then for any distinct $u,v\in \partial H$, we have $\partial i_\gamma^+(u) = \partial i_\gamma^+(v)$ if and only if $(u,v)$ is a leaf of the ending lamination $\Lambda_z$. 
\end{mainthrmC}

\begin{proof}
Suppose first that $\gamma = (z',z)$ is as in \Cref{FormalConventionThroughIdentity} with $1\in \gamma$.
By \Cref{Injectivity}, it suffices to show that if $\partial i_\gamma^+(u) = \partial i_\gamma^+(v)$, then $\lambda = (u,v)\in \Lambda_z$. So, let $u,v\in \partial H$ be distinct points such that 
$\partial i_\gamma^+(u) = \partial i_\gamma^+(v)$. As the set of leaves of $\partial^2 H$ whose endpoints are identified under $\partial i_\gamma^+$ is $H$-invariant, we may assume that $\lambda = (u,v)$ passes through $1\in \Gamma_H$.

Let $\sigma_0: [1,z)\to X(\gamma)^+$ be the isometric lift of $\gamma^+$ into $\Gamma_G$ through the identity as in \Cref{FormalConventionThroughIdentity}. Let $\sigma_e := \Pi_\lambda^{\sigma_0} \cdot \sigma_0$ be the projection of $\sigma_0$ onto $B_{\gamma^+}(\lambda, \sigma_0)$ and set $g_n':= \sigma_e(z_n)$. By \Cref{Analogue of 4.6}, $\sigma_e$ is a $(C, 0)$-quasi-isometric section of $[1,z)$ into $B_{\gamma^+}(\lambda, \sigma_0)$ for some $C \geq 1$. 

By \Cref{Analogue of 4.5}, there exists $C' \geq 1$ such that for any $g\in \sigma_0([1,z))$ and any $[p, q]\subseteq \lambda_{g}$, there exists an extension $[r, q] =: \mu$ of $[p,q]$ in $\lambda_g$ with $d_H(p,r) = 0$ or $1$ and a $(C', 0)$-quasi-isometric section $\sigma$ such that $gr\in \sigma([1, z))$ and $\mu_{r^{-1}g^{-1}g'}$ is a conjugacy minimal representative for all $g'\neq gr$ in $\sigma([1, z))$. Projecting $\sigma$ to $B_{\gamma^+}(\lambda, \sigma_0)$ yields, by \Cref{Analogue of 4.6}, a $(C_2,0)$-quasi-isometric section for some $C_2\geq 1$.

If $\sigma'$ is any $(C_2,0)$-quasi-isometric section, \Cref{Analogue of 4.10} gives that there is some $\alpha > 0$ such that if $\sigma'\subseteq B_{\gamma^+}(\lambda, \sigma_0)$, then there exists some $N \geq 0$ such that for all $n\geq N$, $d_{X(\gamma)^+}(g_n', \sigma'(z_n))\leq \alpha$. Given this $\alpha$, \Cref{Mitra 2.4} guarantees there are some $b> 1$, $A > 0$, and $\eta > 0$ depending on $\alpha$ and $C_2$ such that if $\sigma'([1,z))$ is a $(C_2, 0)$-quasi-isometric section of $[1,z)$ into $X(\gamma)^+$ with $d_{X(\gamma)^+}(\sigma'(z_n), g_n')\geq \eta$, then any path in $i_\gamma^+(\Gamma_H)$ joining $\sigma'(1)$ and $\sigma_e(1)$ has length greater than or equal to $Ab^n$.

Now, let $\lambda^+$ and $\lambda^-$ denote the two closures of the components of $\lambda \setminus \{1\}$. Note that for each $n > 0$, $g_n' \in t_{g_n}\cdot i_\gamma^+(\lambda_{g_n})$. 
Hence, for all $n > 0$ there exists $p_n\in \lambda_{g_n}^-$ and $q_n\in \lambda_{g_n}^+$ such that $d_{X(\gamma)^+}(t_{g_n}\cdot i_\gamma^+(p_n), g_n') = d_{X(\gamma)^+}(g_np_n, g_n') = \eta + 1$ and $d_{X(\gamma)^+}(t_{g_n}\cdot i_\gamma^+(q_n), g_n') = d_{X(\gamma)^+}(g_nq_n, g_n') = \eta$.
By \Cref{Analogue of 4.5}, for each $n>0$ there exists $r_n\in \lambda^-_{g_n}$ with $d_{H}(r_n,p_n) = 0$ or 1 and a $(C',0)$-quasi-isometric section $\sigma_n$ of $[1,z)$ into $X(\gamma)^+$ satisfying the following two conditions: 
\begin{enumerate}
    \item $g_nr_n = \sigma_n(z_n)$
    \item If $\mu^{(n)}$ is the subsegment of $\lambda_{g_n}$ in $\Gamma_H$ joining $r_n$ and $q_n$, then $\mu^{(n)}_{r_n^{-1}g_n^{-1}\sigma_n(z_{m})}$ is a conjugacy minimal representative for all $z_m \neq z_n$.
\end{enumerate}

For each $n>0$, define a new quasi-isometric section $\tau_n(z_i):= t_{g_nq_nr_n^{-1}g_n^{-1}}\cdot \sigma_n(z_i)$ which is obtained by left-translating $\sigma_n$ to go through the point $g_nq_n\in t_{g_n}\cdot i_\gamma^+(\lambda_{g_n})$. We will now project $\sigma_n$ and $\tau_n$ to the set $B_{\gamma^+}(\lambda, \sigma_0)$ to get new quasi-isometric sections which satisfy the hypotheses of \Cref{Analogue of 4.10}. Denote these new $(C_2,0)$-quasi-isometric sections by $\sigma_n':= \Pi_\lambda^{\sigma_0}\cdot \sigma_n$ and $\tau_n':=\Pi_\lambda^{\sigma_0}\cdot \tau_n$.

By \Cref{Analogue of 4.10}, there is some $\alpha$ such that for every index $n>0$, $d_{X(\gamma)^+}(g_k', \sigma_n'(z_k)) \leq \alpha$ as long as $k\geq N$ for some constant $N = N(n)$. So, the $(C_2,0)$-quasigeodesic rays interpolated by $\sigma_n'$ and $\sigma_e$ satisfy the hypotheses of \Cref{Mitra 2.4} since there is some point along these rays where $d_{X(\gamma)^+}(g_k', \sigma_n'(z_k)) \leq \alpha$ and since the rays were defined so that $d_{X(\gamma)^+}(g_n', \sigma_n'(z_n)) = d_{X(\gamma)^+}(g_n',g_nr_n) \geq \eta$. As any path in $i_\gamma^+(\Gamma_H)$ is distance at least $n/C_2$ from any path in $t_{g_n}\cdot i_\gamma^+(\Gamma_H)$, we have that there exists $b > 1$ and $A > 0$ such that the portion of $i_\gamma^+(\lambda)$ between $\sigma'_n(1)$ and $\sigma_e(1) = g_0'$ has length at least $Ab^n$. As the same holds true for the quasigeodesic rays interpolated by $\tau_n'$ and $\sigma_e$, the portion of $i_\gamma^+(\lambda)$ between $\tau'_n(1)$ and $\sigma_e(1) = g_0'$ also has length greater than or equal to $Ab^n$.

Note that for all $n\geq 0$, $\sigma_n(1)$, $\sigma'_n(1)$, $\tau_n(1)$, and $\tau'_n(1)$ all lie in $i_\gamma^+(\Gamma_H)$. Let $[\sigma'_n(1)^*, \tau'_n(1)^*]$ denote the subsegment of $\lambda$ joining $(i_\gamma^+)^{-1} \cdot \sigma'_n(1)$ and $(i_\gamma^+)^{-1}\cdot \tau'_n(1)$. Then, the sequence $\{[\sigma'_n(1)^*, \tau'_n(1)^*]\}$ converges to $\lambda$ in $\widehat{\Gamma_H}$.

Since $d_{X(\gamma)^+}(g_nr_n, g_nq_n)\leq 2\eta + 2$, there exists $\rho > 0$ such that $r_n^{-1}q_n$ is an element of $H$ with $|r_n^{-1}q_n|_H \leq \rho$. Since there are only finitely many of these, we may pass to a subsequence $n_j$ such that $r_{n_j}^{-1}q_{n_j} = h$ where $h$ is some fixed element of $H$. Note that the subsequence $\{[\sigma'_{n_j}(z_0)^*, \tau'_{n_j}(z_0)^*]\}$ also converges to $\lambda$ in $\widehat{\Gamma_H}$.

Let $[\sigma_n(1)^*, \sigma'_n(1)^*]$ denote a geodesic segment in $\Gamma_H$ joinging $(i_\gamma^+)^{-1}\cdot \sigma_n(1)$ and $(i_\gamma^+)^{-1} \cdot \sigma'_n(1)$ and define $[\tau'_n(1)^*, \tau_n(1)^*]$ similarly. Since $\sigma'_n(1)^* = (i_\gamma^+)^{-1}\cdot \Pi_\lambda^{\sigma_0} \cdot i_\gamma^+ (\sigma_n(1))$, we must have that in $\Gamma_H$, $(\sigma_n(1)^*, \tau_n'(1)^*)_{\sigma_n'(1)^*} \leq 2\delta$. Otherwise, there would be a point on $i_\gamma^+(\lambda)$ closer to $\sigma_n(1)$ than $\sigma_n'(1)$, contradicting the definition of $\sigma_n'(1)$ as the projection of $\sigma_n(1)$ to $i_\gamma^+(\lambda)$. For a similar reason, $(\sigma_n'(1)^*, \tau_n(1)^*)_{\tau_n'(1)^*}\leq 2\delta$. So by \Cref{SmallGromovProductsQuasigeodesic}, we have that for all $n$ sufficiently large (so that $d_H(\sigma_n'(1)^*, \tau_n'(1)^*) > 14\delta$), $[\sigma_n(1)^*, \sigma_n'(1)^*]\cup [\sigma_n'(1)^*, \tau_n'(1)^*]\cup [\tau_n'(1)^*, \tau_n(1)^*]$ is a $(1,12\delta)$-quasigeodesic. Thus, for all $n$ sufficiently large, there is some constant $B > 0$ depending only on $\delta$ such that $[\sigma_n(1)^*, \sigma_n'(1)^*]\cup [\sigma_n'(1)^*, \tau_n'(1)^*]\cup  [\tau_n'(1)^*, \tau_n(1)^*]$ lies in a $B$-neighborhood of the geodesic $[\sigma_n(1)^*, \tau_n(1)^* ]$ in $\Gamma_H$.

As the sequence $\{[\sigma_{n_j}'(1)^*, \tau_{n_j}'(1)^*]\}$ converges to $\lambda$, we must also have that the sequence $$\{[\sigma_{n_j}(1)^*,\sigma_{n_j}'(1)^*]\cup [\sigma_{n_j}'(1)^*,\tau_{n_j}'(1)^*]\cup [\tau_{n_j}'(1)^*, \tau_{n_j}(1)^*] \}$$ converges to $\lambda$. In particular, $\{\sigma_{n_j}(1)^*\}$ and $\{\tau_{n_j}(1)^*\}$ must converge to the endpoints of $\lambda$ in $\Gamma_H$. Recall that $\sigma_{n}$ was chosen so that, in particular, $\mu^{(n)}_{r_{n}^{-1}g_n^{-1}\sigma_{n}(1)}$ is a conjugacy minimal representative. Since $\mu^{(n)}_{r_{n}^{-1}g_n^{-1}\sigma_{n}(1)}$ is the label of the geodesic in $\Gamma_H$ between $\sigma_n(1)^*$ and $(g_nq_nr_n^{-1}g_n^{-1}\sigma_n(1))^* = \tau_n(1)^*$, we have that $\{[\sigma_{n_j}(1)^*, \tau_{n_j}(1)^*]\}$ is a sequence of conjugacy minimal representatives of $\phi_{r_{n_j}^{-1}g_{n_j}^{-1}\sigma_{n_j}(1)}(h)$.

Let $\sigma'': [1,z)\to \Gamma_G$ be any quasi-isometric section. Note that for all $n\geq 0$, $\sigma''(z_n)$ and $\sigma_n(1)^{-1}g_nr_n$ are in the same coset of $H$ in $G$. Therefore, $\phi_{r_{n_j}^{-1}g_{n_j}^{-1}\sigma_{n_j}(1)}(h)$ and $\phi_{(\sigma''(z_n))^{-1}}(h)$ have the same conjugacy minimal representatives. Hence, $\lambda = (u,v) \in \Lambda_{z,h}\subseteq \Lambda_z$.

Finally, suppose that $\gamma = (z',z)$ goes through $z_0\in \Gamma_Q$ rather than the identity. Then, $\gamma' := z_0^{-1} \gamma = (z_0^{-1} z', z_0^{-1}z)$ does go through the identity. If $\partial i_\gamma^+(u) = \partial i_\gamma^+(v)$, \Cref{FormalReduction} implies that $\partial i_{\gamma'}^+(\phi_{g_0}(u)) = \partial i_{\gamma'}^+(\phi_{g_0}(v))$. By the above, this implies that $\lambda_{g_0}= (\phi_{g_0}(u), \phi_{g_0}(v)) \in \Lambda_{z_0^{-1}z}$. Finally by \Cref{LambdaZ0}, this implies that $\lambda\in \Lambda_z$ as desired.
\end{proof}

\section{Proof of the main result}\label{MainResult}

We can now prove the main result of the paper, \hyperlink{thrmA}{Theorem A} from the Introduction. Recall that a \textit{dendrite} is a compact, connected, locally connected metrizable space which contains no simple closed curves.

\begin{proposition}[Bowditch \cite{Bow13} c.f. 2.5.2]\label{BowditchDendrite}
Let $\mathcal{X}$ be a bi-infinite hyperbolic stack and let $\mathcal{X}^+$ be the corresponding semi-infinite stack. Then, the Gromov boundary $\partial \mathcal{X}^+$ is a dendrite.
\end{proposition}

If $L\subseteq \partial^2 H$ is an algebraic lamination on $H$, then $\partial H / L$ denotes the quotient space of $\partial H$ by the equivalence relation generated by $L\subseteq \partial^2 H$.

\begin{mainthrmA}\hypertarget{thrmA}{}
Let $1\to H\to G \to Q \to 1$ be a short exact sequence of infinite, finitely generated, word-hyperbolic groups and choose $z\in \partial Q$. Then, the space $\partial H / \Lambda_z$ is homeomorphic to a dendrite.
\end{mainthrmA}

\begin{proof}
Let $\gamma = (z',z)$ be as in \Cref{StackConvention} and let $X(\gamma)^+$ and $i_\gamma^+: \Gamma_H\to X(\gamma)^+$ be as in \Cref{StackInclusionConvention}. Denote by $\pi_z:\partial H\to \partial H/\Lambda_z$ the quotient map. If $a,b\in \partial H$ are such that $\pi_z(a) = \pi_z(b)$, then $\partial i_\gamma^+(a) = \partial i_\gamma^+(b)$ by \Cref{Injectivity}. So, the Cannon-Thurston map $\partial i_\gamma^+: \partial H\to \partial X(\gamma)^+$ quotients through to a map $\tau_z: \partial H / \Lambda_z\to \partial X(\gamma)^+$ with $\partial i_\gamma^+= \tau_z \circ \pi_z$. We will show that $\tau_z$ is a continuous bijection from a compact topological space to a Hausdorff topological space, and thus is a homeomorphism. 

Note that the Gromov boundary of a proper hyperbolic space is compact and metrizable (see for instance \cite{KapBen02}), and so $\partial X(\gamma)^+$ is compact Hausdorff and $\partial H / \Lambda_z$ is compact. As $\partial i_\gamma^+$ is continuous by virtue of being a Cannon-Thurston map (\Cref{i_gamma Cannon-Thurston}) and the quotient map $\pi_z$ is also continuous, the map $\tau_z$ must be continuous. By \hyperlink{thrmB}{Theorem B}, $\partial i_\gamma^+$ is surjective and so $\tau_z$ must also be surjective. If $a',b'\in \partial H / \Lambda_z$ are such that $\tau_z(a') = \tau_z(b') = u\in \partial X(\gamma)^+$, then since $\partial i_\gamma^+$ is surjective, there must exist $a,b\in \partial H$ such that $\partial i_\gamma^+(a) = \partial i_\gamma^+(b) = u$. But by \hyperlink{thrmC}{Theorem C}, this implies that $(a,b)\in \Lambda_z$, and so $\tau_z$ is injective. It now follows that $\tau_z: \partial H / \Lambda_z \to \partial X(\gamma)^+$ is a homeomorphism. Therefore by \Cref{BowditchDendrite}, $\partial H / \Lambda_z$ is a dendrite.
\end{proof}

\footnotesize
\bibliographystyle{plain}
\bibliography{Bib}

\bigskip
\noindent
Department of Mathematics, University of Illinois at Urbana-Champaign\\
1409 W. Green Street, Urbana, IL 61801, U.S.A\\
E-mail: {\tt ecfield2@illinois.edu}

\end{document}